\documentclass[11pt]{amsart}

\usepackage{hyperref}
\usepackage{amssymb,amsmath}
\usepackage[alphabetic,nobysame,abbrev]{amsrefs}
\usepackage{enumerate}
\usepackage{color}
\usepackage{pdfsync}
\RequirePackage{fix-cm}
\usepackage{bm}
\usepackage{mdframed}
\usepackage{cancel}
\usepackage{makecell}
\usepackage{pdfpages}
\usepackage{pgf,tikz}
\usepackage{mathrsfs}
\usetikzlibrary{arrows}
\usetikzlibrary{arrows,decorations.markings}
\usetikzlibrary{matrix}
\usetikzlibrary{backgrounds}

\definecolor{qqqqff}{rgb}{0.,0.,1.}
\definecolor{xdxdff}{rgb}{0.49019607843137253,0.49019607843137253,1.}

\definecolor{qqqqff}{rgb}{0.,0.,1.}

\usepackage{tikz}
\usepackage{graphicx}
\usepackage[cmtip,all]{xy}

\usepackage[draft]{todonotes}
\usepackage[cmtip,all]{xy}

\usepackage{hyperref}
\usepackage[utf8]{inputenc}
\usepackage{amsmath,amsfonts,amssymb,euscript,graphicx,color,tikz}

\theoremstyle{plain}
\newtheorem{theorem}{Theorem}[subsection]
\newtheorem{thm}[theorem]{Theorem}
\newtheorem{lem}[theorem]{Lemma}
\newtheorem{cor}[theorem]{Corollary}
\newtheorem{pro}[theorem]{Proposition}

\theoremstyle{definition}
\newtheorem{DEF}[theorem]{Definition}

\newtheorem{exa}[theorem]{Example}

\newtheorem{rem}[theorem]{Remark}

\newtheorem{parag}[theorem]{{}}

\newcommand\pref[1]{\textbf{\ref{#1}}}
\numberwithin{equation}{section}

\newcommand{\sub}{\subseteq}

\newcommand{\la}{Lie algebra }

\newcommand{\fm}{(\cdot,\cdot)}
\newcommand{\fh}{\mathfrak{h}}

\newcommand{\K}{\mathbb{K} }

\newcommand{\Z}{\mathbb{Z} }

\newcommand{\fg}{\mathfrak{g}}

\newcommand{\ep}{\hfill$\Box$}

\def\ad{\hbox{ad}}
\def\andd{\quad\hbox{and}\quad}
\def\sg{\sigma}
\def\a{\alpha}
\def\b{\beta}
\def\lam{\lambda}
\def\Lam{\Lambda}
\def\ep{\epsilon}
\def\andd{\quad\hbox{and}\quad}
\def\supp{\hbox{supp}}
\def\id{\hbox{id}}
\def\Aut{\hbox{Aut}}

\def\End{\mathrm{End}}

\def\g{\mathfrak{g}}

\def\mod{\hbox{mod}}
\def\andd{\quad\hbox{and}\quad}
\def\ind{\hbox{ind}}

\def\v{{\mathcal V}}
\def\u{{\mathcal U}}

\def\fm{(\cdot,\cdot)}
\def\a{\alpha}

\def\sub{\subseteq}

\def\lam{\lambda}
\def\Lam{\Lambda}

\def\1k{\frac{1}{k}}

\def\la{\langle}
\def\ra{\rangle}

\def\d{\delta}

\def\b{\beta}

\def\qed{\hfill$\Box$}

\def\sg{\sigma}

\def\hh{{\mathcal H}}

\def\sg{\sigma}

\usepackage{verbatim} 

\def\quadd{\quad\quad}

\def\ad{\hbox{ad}}

\def\be{{\bf e}}
\def\bq{{\bf q}}

\def\bbbz{{\mathbb Z}}
\def\Z{{\mathbb Z}}

\def\bbbr{{\mathbb R}}

\def\bbbk{{\mathbb K}}

\def\aa{\mathcal A}

\def\tr{\hbox{tr}}
\def\cc{{\mathcal C}}
\def\dd{\mathcal D}
\def\ep{\epsilon}

\def\ll{{\mathcal G }}

\def\jj{{\mathcal J}}
\def\ll{\mathcal L}

\def\gl{\frak{g}\ell}

\def\supp{\hbox{supp}}

\def\bb{{\mathcal B}}
\def\uu{{\mathcal U}}

\def\proof{{\noindent\bf Proof. }}

\def\g{\mathfrak{g}}

\def\span{\hbox{span}}

\def\cent{\mathrm{Cent}}

\def\ev{\hbox{\bf ev}}
\def\scd{\mathrm{SCDer}}

\def\DynkinNodeSize{1.5mm}
\def\DynkinArrowLength{2mm}
\tikzset{
	dnode/.style={
		circle,
		inner sep=0pt,
		minimum size=\DynkinNodeSize,
		fill=white,
		draw},
	middlearrow/.style={
		decoration={markings,
			mark=at position 0.8 with
			{\draw (0:0mm) -- +(+140:\DynkinArrowLength); \draw (0:0mm) -- +(-140:\DynkinArrowLength);},
		},
		postaction={decorate}
	},
	leftrightarrow/.style={
		decoration={markings,
			mark=at position 0.999 with
			{
				\draw (0:0mm) -- +(+135:\DynkinArrowLength); \draw (0:0mm) -- +(-135:\DynkinArrowLength);
			},
			mark=at position 0.001 with
			{
				\draw (0:0mm) -- +(+45:\DynkinArrowLength); \draw (0:0mm) -- +(-45:\DynkinArrowLength);
			},
		},
		postaction={decorate}
	},
	sedge/.style={
	},
	dedge/.style={
		middlearrow,
		double distance=0.6mm,
	},
	tedge/.style={
		middlearrow,
		double distance=1.0mm+\pgflinewidth,
		postaction={draw}, 
	},
	infedge/.style={
		leftrightarrow,
		double distance=0.5mm,
	},
}

\begin{document}

%
%
\title{Chevalley involutions for Lie tori and extended affine Lie algebras}

\author{Saeid Azam, Mehdi Izadi Farhadi}
\address
{Department of Pure Mathematics\\Faculty of Mathematics and Statistics\\
	University of Isfahan\\ P.O.Box: 81746-73441\\ Isfahan, Iran, and\\
	School of Mathematics, Institute for
	Research in Fundamental Sciences (IPM), P.O. Box: 19395-5746.} \email{azam@ipm.ir, azam@sci.ui.ac.ir}
\address
{Department of Pure Mathematics\\Faculty of Mathematics and Statistics\\
	University of Isfahan\\ P.O.Box: 81746-73441\\ Isfahan, Iran.} \email{m.izadi@ipm.ir}
\thanks{This work is based upon research funded by Iran National Science Foundation (INSF) under project No.  4001480.}
\thanks{This research was in part carried out in
	IPM-Isfahan Branch.}


\begin{abstract}
In finite-dimensional simple Lie algebras and affine Kac-Moody Lie algebras, Chevalley involutions are crucial ingredients of the modular theory.
Towards establishing the modular theory for extended affine Lie algebras, we investigate the existence of ``Chevalley involutions" for Lie tori and extended affine Lie algebras.  
We first discuss how to lift a Chevalley involution from the centerless core
which is characterized to be a centerless Lie torus to the core and then to the entire extended affine Lie algebra. {We then prove by a  type-dependent argument the existence of Chevalley involutions for centerless Lie tori.}
\end{abstract}
\maketitle

\setcounter{equation}{-1}
\section{\bf Introduction}\setcounter{equation}{0}
\label{intro}
\markboth{S. Azam, M. Izadi Farhadi}{Chevalley Involutions}

The main goal of this study is to investigate the existence of Chevalley involutions for extended affine Lie algebras (EALAs).
Extended affine Lie algebras are a class of mostly infinite dimensional Lie algebras that generalize the class of finite and affine Lie algebras. 
This class has been under intensive investigation for the past few decades.
 For the basic
structure theory of extended affine Lie algebras, we refer the interested reader to \cite{AABGP97} and \cite{Neh11}.

This work is a continuation of \cite{AFI22} in which
the concept of modular theory for extended affine Lie algebras was initiated. By a modular theory, we mean a higher nullity version of the seminal procedure in finite dimensional complex simple Lie algebras due to C. Chevalley which provides a framework for passing from complex numbers to any field of positive characteristic, leading to the Chevalley groups. The main ingredients of this procedure are Chevalley involutions and Chevalley bases. Roughly speaking, by means of a Chevalley involution which acts on the Cartan subalgebra as  minus the identity, one constructs a Chevalley basis whose $\bbbz$-span provides a $\bbbz$-form {(integral form)} for the ground complex Lie algebra. Then by tensoring the provided $\bbbz$-form with any field of positive characteristic one obtains a modular Lie algebra
that is independent of the choice of the Chevalley base. One then proceeds to get an integral form for the corresponding universal enveloping algebra which is of great importance in the representation theory of Chevalley groups, see \cite{Kos66} and \cite{Hum72}.

The existence of Chevalley involutions is a focal point in modular theory. Naively, such involutions permit the lift of arguments that are valid locally in certain subalgebras to arguments that are valid globally. In \cite{AFI22}, the authors construct a $\bbbz$-form for the core of the considered extended affine Lie algebra, by assuming that it admits a Chevalley involution. In finite and affine theory, the existence of such involutions is guaranteed by Serre-type presentations. However, the situation for general extended affine Lie algebras is not clear at all. In this work, we investigate the existence of Chevalley involutions for (tame) extended affine Lie algebras.

H. Garland extended the theory of integral forms to untwisted affine Lie algebras by providing Chevalley bases for these algebras and presenting integral forms of corresponding enveloping algebras along with explicit $\bbbz$-bases of them. Also he constructed Chevalley groups associated to standard representations of untwisted affine Lie algebras, see \cite{Gar78}, \cite{Gar80}. In \cite{Mit85}, D. Mitzman gives, by a type-dependent approach, a $\bbbz$-basis for the universal enveloping algebra corresponding to a twisted affine Lie algebra. The construction of $\bbbz$-forms for the class of all Kac-Moody Lie algebras is due to J. Tits \cite{Tit81}, \cite{Tit82}.

We recall that an extended affine Lie algebra is a triple $(E,\fm,\hh)$ consisting of a
complex Lie algebra $E$, a Cartan subalgebra $\hh$ of $E$, and a symmetric invariant non-degenerate bilinear form $\fm$ on $E$, satisfying certain natural axioms, see Definition \ref{EALA}. Let $R$ be the set of roots of $E$ with respect to $\hh$. Then we have $R=R^\times\cup R^0$ where $R^\times$ is the set of ``non-isotropic'' roots of $E$ and $R^0$ is the set of ``isotropic roots'' of $E$. 
To explain our strategy in investigating Chevalley involutions for $E$, we need to first recall a few terms.
The ``core'' $E_c$ of $E$ is by definition the subalgebra of $E$ generated by non-isotropic root spaces $E_\a$, $\a\in R^\times$.
The centerless core $E_{cc}$ is the quotient Lie algebra $E_c/Z(E_c)$, where $Z(E_c)$ is the center of $E_c$.
  It is very well understood that the centerless core plays a central role in the classification and realization of extended affine Lie algebras as well as other parts of the theory, see \cite{AG01}. The centerless core is characterized as a ``centerless Lie torus'', see Definition \ref{def1}. Centerless Lie tori are defined by Y. Yoshii \cite{Yos06}. 
  In brief, the centerless core of an extended affine Lie algebra is a centerless Lie torus. Conversely, starting from a centerless Lie torus $\ll$ there is a prescribed way due to E. Neher \cite{Neh04} to construct an extended affine Lie algebra whose centerless core is $\ll$. We explain this construction briefly here.

Let $\ll$ be a centerless Lie torus equipped with a symmetric invariant non-degenerate form $\fm$; see \cite[Theorem 5.2]{Yos06}. Let $D$ be a ``permissible subalgebra" of skew centroidal derivations $\scd(\ll)$ of $\ll$, and ${D^{gr}}^\star$ be its graded dual, viewed as a $D$-module via the contragredient action, see \pref{parag2} and Definition \ref{permis} for details. Let $\kappa$ be an {``affine cocycle''} on $D$, see \pref{parag1}. Consider the vector space $E=E(\ll,D,\kappa)=\ll\oplus{D^{gr}}^\star\oplus D$.
Then one endows $E$ with a natural Lie bracket induced by the Lie bracket on $\ll$, the commutator bracket on $D$, the action of $D$ on ${D^{gr}}^\star$, the natural action of $D$ on $\ll$ as derivations, the affine cocycle $\kappa$, and an involved $2$-cocycle invoked by $D$ and $\fm$, see \pref{parag3}. One also endows $E$ with a bilinear form $\fm$, induced by the form on $\ll$ and the dual paring of $D$ and ${D^{gr}}^\star.$ set $\hh=\ll^0_0\oplus{D^{0}}^\star\oplus D^0$.
According to \cite[Theorem 16]{Neh04}, $(E,\fm,\hh)$ is an extended affine Lie algebra with $E_c=\ll\oplus{D^{gr}}^\star$, conversely any extended affine Lie algebra is graded-isomorphic to $E(\ll,D,\kappa)$ for some centerless Lie torus $\ll$, some permissible subalgebra $D$, and an affine cocycle $\kappa$. 

We now return to the question of the existence of Chevalley involutions for an extended affine Lie algebra $E(\ll,D,\kappa)$. Suppose first that $\ll$ is equipped with such an involution $\tau$, see Definition \ref{1def}. We extend $\tau$ to an involution $\bar\tau$ for the vector space $\ll\oplus\End(\ll)^\star\oplus\End(\ll)$ in a natural manner. Then $\bar\tau$ restricts to an involution $\bar\tau$ of $\ll\oplus{\scd(\ll)^{gr}}^\star\oplus\scd(\ll)$ as a Lie algebra. Now the question is if $\bar\tau$ preserves the permissible subalgebra $D$ and its graded dual. The answer in general is negative, see Example \ref{sky-10}(ii). To get a better picture of the situation, we set $D_\tau:=\bar\tau(D)$ and we define the bilinear map
$\kappa_\tau:D_\tau\times D_\tau\rightarrow {({D_\tau^0})^\star}$ by $\kappa_\tau=\bar\tau\circ
\kappa\circ(\bar\tau\times\bar\tau).$
It turns out that $D_\tau$ is a permissible subalgebra of $\scd(\ll)$ and $\kappa_\tau$ is an affine cocycle on $D_\tau$, see Lemma \ref{1lem}. Therefore the Chevalley involution $\tau$ for $\ll$ leads us to a Lie algebra isomorphism $\bar\tau$ from
$E(\ll,D,\kappa)$ onto $E^\tau:=E(\ll,D_\tau,\kappa_\tau)$ which maps the core of $E$ onto the core of $E^\tau$. Moreover, if $D$ is $\tau$-invariant, namely $D_\tau=D$, then $\bar\tau$ restricts to a Chevalley involution of $E_c$.
Furthermore, if the pair $(D,\kappa)$ is $\tau$-invariant, meaning $(D_\tau,\kappa_\tau)=(D,\kappa)$, then
$\bar\tau$ restricts to a Chevalley involution of $E$, see Theorem \ref{induce}. In Lemma \ref{sky19}, we give a necessary and sufficient condition for a permissible subalgebra $D$ to be $\tau$-invariant. It turns out that the special cases $D=D^0$ or $D=\scd(\ll)$ are both $\tau$-invariant, and so the existence of involution $\bar\tau$ on {$E_c$} is guaranteed.
As far as we know, these two cases cover all specific examples of extended affine Lie algebras {appeared} in the literature, see remark \ref{remneg}, except the ones we construct in Example \ref{sky-10}. The question; ``Which affine cocycles are invariant?" is more subtle since we don't have enough examples and our information about affine cocycles is limited. 
Fortunately, the known {examples are} $\tau$-invariant with respect to any Chevalley involution $\tau$ for $\ll$, see Remark \ref{remneg}(ii).

{We now discuss the existence of Chevalley involutions for centerless Lie tori of reduced types. Fortunately, due to the works of \cite{BGK96},\cite{BGKN95}, \cite{Yos00} and \cite{AG01} the structures and coordinate algebras of centerless Lie tori of reduced types are very well understood. 
	In fact a centerless Lie torus $\ll$ is coordinatized by a unital algebra that is a Jordan algebra, an alternative algebra or an associative algebra, depending on its type $\Delta$, see Section \ref{tori}. In simply laced cases, types $A_\ell(\ell\geq3)$, $A_2$ and $A_1$ can be constructed using $\frak{sl}_{\ell+1}$, $\frak{psl}_3$ and Tits-Kantor-Koecher (TKK) construction, respectively, and types $D,E$ are untwisted centerless Lie tori. In non-simply laced cases, types $B_\ell(\ell\geq3)$, $F_4$ and $G_2$ can be constructed using generalized Tits construction and type $C_\ell (\ell\geq2)$ can be constructed using Tits-Kantor-Koecher construction, see Section \ref{Chevalley involutions for centerless Lie tori}. We, therefore, start by showing that the coordinate algebra of $\ll$ admits a pre-Chevalley involution, meaning that it admits an involution that reverses the grading, see Section \ref{tori}. We then, {by a type-dependent procedure}, show that this involution can be lifted to a Chevalley involution for the corresponding tensor structure, $\frak{sl}_{\ell+1}$ structure, $\frak{psl}_3$ structure,
	TKK structure, or generalized Tits structure, depending on the type. 
		{This completes the problem of the existence of Chevalley involutions for the considered centerless Lie tori, see Theorems \ref{thm-2} and \ref{thm-3}.} The results mentioned above covers the content of Sections {2-5}.}

In Section \ref{sec7}, we discuss the existence of Chevalley involutions for Lie tori that are finitely generated over their centroid (fgc Lie tori).
Any centerless Lie torus of a type different from $A_\ell$ is fgc, and any such algebra is a multi-loop algebra, see Section \ref{sec7} for details and definitions. In Subsection \ref{subsec7}, the existence of Chevalley involutions for multi-loop algebras is considered.  

{We conclude the introduction by mentioning that the study of involutions is closely related to the study of compact forms of involutive Lie algebras. In \cite{Gao96} a uniform approach is considered for the study of involutive Lie algebras graded by finite root systems of simply laced types.}

\section{\bf Preliminaries}\setcounter{equation}{0}
All vector spaces and algebras in this work will be  over the field $\bbbk$ of complex numbers, except when indicated otherwise. If $S$ is a subset of a group, we denote by $\langle S\rangle$ the subgroup generated by $S$. Also throughout this work $\aa_{[n]}=\bbbk[x_1^{\pm1},\ldots,x_n^{\pm1}]$ denotes the algebra of Laurent polynomials over $\K$ in $n$-variables.

\subsection{Extended affine Lie algebras}
We begin by recalling the definition of an ``extended affine Lie algebra". Since in this work, we are interested in ``tame" extended affine Lie algebras, the tameness condition is included in the definition. For details, we refer the reader to \cite{AABGP97} and \cite{Neh11}.
\begin{DEF}\label{EALA}
	An {\it extended affine Lie algebra} is a triple
	$(E,\fm,\hh)$ consisting of a Lie algebra $E$,
	a non-trivial subalgebra $\hh$ and a bilinear form $\fm:E\times E\rightarrow\bbbk$ satisfying the following six axioms:
	
	(A1) $\fm$ is invariant and non-degenerate on $E$,
	
	(A2) $\hh$ is finite-dimensional and we have the root space decomposition 
	$E=\sum_{\a\in\hh^\star}E_\a$ with
	$E_\a=\{x\in E\mid [h,x]=\a(h)x\hbox{ for all }h\in\hh\},$
	and $E_0=\hh.$ 
	
To give the next axioms, we denote the root system of $E$ with respect to $\hh$ by $R$, namely  $R=\{\a\in\hh^\star\mid E_\a\neq\{0\}\}$.
From (A1)-(A2), it follows that the form $\fm$ on $E$ restricted to $\hh$ is non-degenerate and so it transfers
	to $\hh^\star$ by $(\a,\b):=(t_\a,t_\b)$ where $t_\a\in\hh$
	is the unique element satisfying $\a(h)=(h,t_\a)$, $h\in\hh$. Then a root $\a$ is called {\it isotropic} if
$(\a,\a)=0$ and {\it non-isotropic} otherwise.
	We denote by $R^0$ and $R^\times$ the set of isotropic and non-isotropic roots of $R$, respectively. The subalgebra $E_c$ of $E$ generated by non-isotropic root spaces is called the {\it core} of $E$. 
	
	(A3) For $\a\in R^\times$ and $x\in E_\a$, $\ad(x)$ is locally nilpotent on $E$.
	
	(A4) $E$ is {\it tame}, namely the centralizer of $E_c$ in $E$ is contained in $E_c$.
	
	(A5) The $\bbbz$-span of $R$ in $\hh^\star$ is a free abelian group of finite rank.
	
	(A6) $R^\times$ is indecomposable in the sense that if $R^\times=R_1\cup R_2$ with $(R_1,R_2)=0$, then $R_1=\emptyset$ or $R_2=\emptyset.$
%
	 \end{DEF}
 
\begin{parag}\label{p0}
	Let $(E,\fm,\hh)$, or simply $E$, be an extended affine Lie algebra with root system $R$. Set $\v=\span_{\bbbr}R$ and {$\v^0=\span_{\bbbr}R^0$}. Let $\bar{\;}:\v\rightarrow\bar\v:=\v/\v^0$.
It turns out that the image $\bar R$ of $R$ in $\bar\v$ is an irreducible finite root system in $\bar\v$. One can find an isomorphic preimage $\dot R$ of $\bar R$ in $\v$, under $\bar{\;}$, such that $\dot R$ is an irreducible finite root system in its real span. Then $R\sub\dot R+\Lam$, where $\Lam=\la R^0\ra$ is a free abelian group of rank equal to $\dim\v^0$.  The {\it type} and the {\it rank} of $R$ (or $E$) is by definition the type and the rank of $\bar R$, respectively.
\end{parag} 

 The quotient Lie algebra
$E_{cc}:=E_c/Z(E_c)$ is called the {\it centerless core} of $E$, where here $Z(E_c)$ denotes the center of $E_c$. Most structural properties of an extended affine Lie algebra is encoded in its centerless core. Centerless core of an extended affine Lie algebra is characterized as a ``Lie torus".

\subsection{Lie tori}\label{Lie tori}
Lie tori were introduced by Y. Yoshii
in \cite{Yos06} to give a characterization of extended affine Lie algebras and were further studied by E. Neher in \cite{Neh04-1}. Here we recall the definition of a Lie torus.

Let $\Lam$ be a free abelian group of finite rank and let $\Delta$ be a finite irreducible root
system with root lattice $Q=\span_\bbbz(\Delta)$. For $\a,\b\in\Delta,\a^\vee$ is the
coroot of $\a$, and $\langle\b,\a^\vee\rangle$ is the corresponding  Cartan integer. We denote by $\Delta_{\ind}=\{0\}\cup\{\a\in\Delta^\times\mid(\a/2)\notin\Delta\}$ the subsystem of ``indivisible'' roots of
$\Delta$.

\begin{DEF}\label{def1}
	A Lie $\Lam$-torus of type $\Delta$ is a Lie algebra $\ll$ over $\bbbk$ satisfying the following conditions (LT1)-(LT4):\\
	(\textbf{LT1}) $\ll$ is a $Q\times\Lam$-graded Lie algebra;
	$$\ll=\bigoplus_{(\a,\lam)\in Q\times\Lam} \ll_\a^\lam\quad\text{with}\quad\ll_\a^\lam=0 \;\text{if} \;\a\notin\Delta,$$\\
	and $$[\ll_\a^\lam,\ll_\b^\mu]\subset \ll_{\a+\b}^{\lam+\mu}.$$
	(\textbf{LT2}) For $\a\in\Delta^\times$ and $\lam\in\Lam$ we have
	\begin{itemize}
		\item [(i)] $\dim \ll_\a^\lam\leq1$, with $\dim \ll_\a^0=1$ if $\a\in\Delta_{\ind}$,
		\item[(ii)] if $\dim \ll_\a^\lam=1$ then there exist elements $e_\a^\lam\in \ll_\a^\lam$ and $f_\a^\lam\in \ll_{-\a}^{-\lam}$ such that
		$$\ll_\a^\lam=\bbbk e_\a^\lam,\quad\ll_{-\a}^{-\lam}=\bbbk f_\a^\lam,$$
		and{ 
		$$[[e_\a^\lam,f_\a^\lam],x^\mu_\b]=\langle\b,\a^\vee\rangle x_\b^\mu$$
		for $\b\in\Delta,$ $\mu\in\Lam$, $x_\b^\mu\in \ll_\b^\mu.$}
		\end{itemize}
	(\textbf{LT3}) For $\lam\in\Lam$ we have $\ll_0^\lam=\sum_{\a\in\Delta^\times,\mu\in\Lam} [\ll_\a^\mu,\ll_{-\a}^{\lam-\mu}]$.\\
	(\textbf{LT4}) $\Lam=\langle\supp_\Lam(\ll)\rangle$, where $\supp_\Lam(\ll)=\{\lam\in\Lam\mid \ll_\a^\lam\neq 0\;\text{for some} \;\a\in\Delta\}$.
	
{The Lie torus $\ll$ is called {\it centreless} if $\ll$ has trivial center, and is called {\it invariant}, if $\ll$ has an invariant non-degenerate symmetric
bilinear form $\fm$ which is graded, meaning that
$$ (\ll^\lam_\a,\ll^\mu_\b)=0\quad\hbox{unless}\quad
(\a,\lam)=(-\b,-\mu).
$$}
\end{DEF}
The rank of $\Lam$ is called the \textit{nullity} of $\ll$. When $\Delta$ and $\Lam$ are fixed, we simply  say that $\ll$ is a Lie torus. For $\lam\in\Lam$ and $\a\in Q$, we set 
$$\ll^\lam:=\bigoplus_{\a\in Q} \ll_\a^\lam\andd\ll_\a:=\bigoplus_{\lam\in\Lam} \ll_\a^\lam,$$
then $\ll=\bigoplus_{\lam\in\Lam}\ll^\lam$ is $\Lam$-graded and $\ll=\bigoplus_{\a\in Q}\ll_\a$ is $Q$-graded.

\begin{parag}\label{parag4}
	{Assume that $E$ is an extended affine Lie algebra with root system $R$. Consider the finite root system $\dot R$, and the free abelian group $\Lam$ as in \pref{p0}. We denote the type of $\dot R$ by $\Delta$. The following result associates a (centreless) Lie torus to $E$.}
\end{parag}

\begin{pro}\cite[Proposition 12]{Neh04}\cite[Proposition 1.28]{AG01}\label{pro4}
	The core of $E$ is a Lie $\Lam$-torus of type $\Delta$ with $Q\times\Lam$-grading $E_c=\bigoplus_{\a\in\dot R,\lam\in\Lam}(E_c)_\a^\lam$, where $(E_c)_\a^\lam=E_c\cap E_{\a+\lam}$. {In particular $E_{cc}=E_c/Z(E_c)=\bigoplus_{\a\in\dot R,\lam\in\Lam}(E_{cc})_\a^\lam$ is a centreless Lie $\Lam$-torus of type $\Delta$ with $(E_{cc})_\a^\lam=\pi((E_c)_\a^\lam)$, where $\pi:E_c\rightarrow E_{cc}$ is the canonical projection defined by $\pi(x)=x+Z(E_c)$.}
\end{pro} 

Certain subalgebras of derivation algebra of a Lie torus play a crucial role in the sequel.

\subsection{Centroidal and skew centroidal derivations}\label{Centroidal and skew centroidal derivations}

{We describe a construction that associates to a centerless Lie torus $\ll$  a class of extended affine Lie algebras $E(\ll,\dd,\kappa)$, where $D$ is a graded subalgebra of derivation algebra of $\ll$ and $\kappa$ is a so called an ``affine cocycle''. As a vector space, we have $E(\ll,D,\kappa)=\ll\oplus{D^{gr}}^\star\oplus D$, where ${D^{gr}}^\star$ is the graded dual of $D$. Below, we will describe the details of this construction.}

\begin{DEF}
	The \textit{centroid} of a Lie algebra $\ll$ over any commutative ring $R$, denoted $\cent_R(\ll)$, is the set of all $R$-linear endomorphisms of $\ll$ that commute with left and right multiplication by elements of $\ll$, i.e.
	$$\cent_R(\ll)=\{\chi\in\mathrm{End}_R(\ll)\mid\chi([x,y])=[\chi(x),y]=[x,\chi(y)],\forall x,y\in \ll\}.$$
Then	$\ll$ can be naturally viewed as a module over $\cent_R(\ll)$ via $\chi.x=\chi(x)$ for all $x\in\ll,\chi\in\cent_R(\ll)$.
\end{DEF}

\begin{DEF}\label{fgc}
	An EALA is called \textit{fgc} if its centerless core is fgc, i.e. the centerless core is finitely generated as a module over its centroid.
\end{DEF}

\begin{parag}
 Let $\ll$ be a Lie $\Lam$-torus of type $\Delta$ over $\bbbk$ of nullity $n$. Then its centroid $\cent_\bbbk(\ll)$ is also $\Lam$-graded;
$$\cent_\bbbk(\ll)=\bigoplus_{\lam\in\Lam}\cent_\bbbk(\ll)^\lam\quad\text{with}\quad\dim_\bbbk\cent_\bbbk(\ll)^\lam\leq 1,$$
where $\cent_\bbbk(\ll)^\lam=\{\chi\in\cent_\bbbk(\ll)\mid\chi(\ll^\mu)\subseteq \ll^{\lam+\mu},\;\text{for all}\;\mu\in\Lam\}$. Now let $\Gamma=\{\lam\in\Lam\mid \cent_\bbbk(\ll)^\lam\neq 0\}$. The following result from \cite{Neh04} justifies to call $\Gamma$ the \textit{centroidal grading group} of $\ll$.
\end{parag}

\begin{thm}\label{thm9}
	Let $\ll$ be a centerless Lie $\Lam$-torus of type
$\Delta$. Then
	\begin{itemize}
		\item [(a)]  $\Gamma$ is a subgroup of $\Lam$, and $\cent_\bbbk(\ll)$ is isomorphic to the group ring $\bbbk[\Gamma]$, hence
		to a Laurent polynomial ring in finitely many variables.

		\item[(b)] $\ll$ is a free $\cent_\bbbk(\ll)$-module, and if $\Delta\neq A_\ell$, then $\ll$ has finite rank as
		a $\cent_\bbbk(\ll)$-module.
	\end{itemize}
\end{thm}

\begin{parag}
Let $\ll$ be a centerless Lie $\Lam$-torus of type $\Delta$ of nullity $n$, with centroidal grading group $\Gamma$. Thus we can write
$$\cent_\bbbk(\ll)=\bigoplus_{\mu\in\Gamma}\bbbk\chi^\mu,$$
where $\chi^\mu$ acts on $\ll$ as an endomorphism of degree $\mu$ and $\chi^\mu\chi^\nu=\chi^{\mu+\nu}$.
For any $\theta\in\text{Hom}_\bbbz(\Lam,\bbbk)$, define derivation $\partial_\theta$ of $\ll$ by
$$\partial_\theta(x^\lam)=\theta(\lam)x^\lam\;\;\text{for}\;\lam\in\Lam,x^\lam\in \ll^\lam.$$
The derivation $\partial_\theta$ is called a \textit{degree derivation} of $\ll$. Put $$\mathcal{D}:=\{\partial_\theta\mid\theta\in\text{Hom}_\bbbz(\Lam,\bbbk)\},$$ the set of all degree derivations. Denote by
$$\mathrm{CDer}_\bbbk(\ll):=\cent_\bbbk(\ll)\mathcal{D}=\bigoplus_{\mu\in\Gamma}\chi^\mu\mathcal{D},$$
called the \textit{centroidal derivations} of $\ll$. It follows that $\mathrm{CDer}_\bbbk(\ll)$ is a $\Gamma$-graded subalgebra of the derivation algebra
$\mathrm{Der}_\bbbk(\ll)$ of $\ll$ with
\begin{equation}\label{eq7}
	[\chi^\mu\partial_\theta,\chi^\nu\partial_\psi]=\chi^{\mu+\nu}(\theta(\nu)\partial_\psi-\psi(\mu)\partial_\theta).
\end{equation}
\end{parag}

\begin{parag}\label{parag2}
Fix a non-degenerate invariant $\Lam$-graded bilinear form $\fm_\ll$ on $\ll$ (see \cite[Theorem 5.2]{Yos06}). A derivation $d$ satisfying $(d(x),x)_\ll=0$ for all $x\in \ll$, {or equivalently $(d(x),y)=-(x,d(y))$ for all $x,y\in\ll$,} is called a \textit{skew derivation} of $\ll$. Set
\begin{eqnarray*}
	\mathrm{SCDer}_\bbbk(\ll)&:=&\{d\in\mathrm{CDer}_\bbbk(\ll)\mid(d(x),x)_L=0\;\text{for all}\;x\in \ll\}\\
	&=&\bigoplus_{\mu\in\Gamma}\mathrm{SCDer}_\bbbk(\ll)^\mu=\bigoplus_{\mu\in\Gamma}\chi^\mu\{\partial_\theta\in\dd\mid\theta(\mu)=0\}.
\end{eqnarray*}
Then $\mathrm{SCDer}_\bbbk(\ll)$ is a $\Gamma$-graded subalgebra of $\mathrm{CDer}_\bbbk(\ll)$, called the {\it skew centroidal derivations} of $\ll$. Note that $\mathrm{SCDer}_\bbbk(\ll)^0=\mathcal{D}.$

For a graded subalgebra $D=\sum_{\mu\in\Gamma}D^\mu$ of
$\mathrm{SCDer}_\bbbk(\ll)$, we denote its graded dual with
$ D^{gr^\star}=\sum_{\mu\in\Gamma}(D^\mu)^\star$ with grading
$(D^{gr^\star})^\mu=(D^{-\mu})^\star$, and we consider it as a $D$-module by the contragredient action, i.e. 
 $$(d.\varphi)(d^{'})=\varphi([d^{'},d])\:\text{for}\;d^{'},d\in D,\varphi\in D^{gr*},$$
	where $\varphi\in(D^\mu)^\star$ is viewed as a linear form on $D$ by $\varphi|_{D^\nu}=0$ for $\nu\neq\mu$. 
\end{parag}

\subsection{Construction of extended affine Lie algebras}\label{Construction of EALAS}
We are ready now to give a general construction of extended affine Lie algebras starting from a centerless Lie torus, see \cite{Neh11}. We begin with the following definition.
\begin{parag}\label{parag1}
Let  $D=\bigoplus_{\mu\in\Gamma}D^\mu$ be a $\Gamma$-graded subalgebra of  $\mathrm{SCDer}_\bbbk(\ll)$, and  $\kappa:D\times D\rightarrow D^{gr*}$ be a bilinear map. Assume that $D$ satisfies the following two conditions:

(i)  the canonical evaluation map $\text{ev}:\Lam\rightarrow(D^0)^\star$ defined by
	$$\text{ev}(\lam)(\partial_\theta)=\theta(\lam),\;\lam\in\Lam,$$
	is injective {and has discrete image}. Note that, we conclude from (i) that $D^0$ is an ad-diagonalizable subalgebra of $D$ with weight spaces $D^\mu$, $\mu\in\Gamma$.
We also note that $D^0\subseteq\mathrm{SCDer}_\bbbk(\ll)^0=\mathcal{D}$.
	
	(ii) $\kappa$ is an ``affine cocycle''   which is graded and invariant. This means that,  $\kappa$ is a bilinear map satisfying
$$\begin{array}{c}
\kappa(d,d)=0,\;\;{\sum_{(i,j,k)\circlearrowleft}\kappa([d_i,d_j],d_k)=\sum_{(i,j,k)\circlearrowleft}d_i\cdot\kappa (d_j,d_k)},\vspace{2mm}\\
\kappa(D^{\mu_1},D^{\mu_2})\subseteq(D^{-\mu_1-\mu_2})^\star\;\;\text{and}\;\;\kappa(d_1,d_2)(d_3)=\kappa(d_2,d_3)(d_1),\vspace{2mm}\\
	\kappa(D^0,D)=0,
	\end{array}
	$$
	for $d,d_1,d_2,d_3\in\ D$.
	 {Note that by $(i,j,k)\circlearrowleft$, we mean $(i,j,k)$ is a cyclic permutation of $(1,2,3)$}.
\end{parag}

\begin{DEF}\label{permis}
{We call a subalgebra $D$ of $\scd(\ll)$  {\it permissible}, if it satisfies \pref{parag1}(i).}
\end{DEF}

\begin{rem}
{Extended affine Lie algebras are defined originally over the field of complex numbers.  This notion, later was generalized by E. Neher to fields of characteristic 0, see \cite[Definition 11]{Neh04}.  In this general case, the canonical evaluation map doesn't need to have discrete image.}
\end{rem}

\begin{exa}\label{eeax1}
	(i) Let $(\u,\u^\prime,\Gamma^\prime)$ be a triple consisting of $\bbbk$-subspaces $\u,\u^\prime$ of $\mathrm{Hom}_\bbbz(\Lam,\bbbk)$ and a subgroup $\Gamma^\prime$ of $\Gamma$.  We set
	$${(\u,\u',\Gamma')_\mu}:=\left\{\begin{array}{ll}
		\u&\hbox{if }\mu=0,\\
		\{\theta\in \u^\prime\mid\theta(\mu)=0\} &\hbox{if }\mu\in\Gamma'\setminus\{0\},\\
		\{0\}&\hbox{otherwise,}
		\end{array}\right.
	$$
	and
	$$D^\mu:=\chi^\mu\{\partial_\theta\mid\theta\in{(\u,\u',\Gamma')_\mu}\},\qquad(\mu\in\Gamma).$$
Then one can see from (\ref{eq7}) that
	$$D_{(\u,\u^\prime,\Gamma^\prime)}:=D^0\oplus\sum_{0\neq\mu\in\Gamma^\prime}D^\mu$$
	is a $\Gamma$-graded subalgebra of $\mathrm{SCDer}_\bbbk(\ll)$. 

	(ii) Note that if $\u$ {in part (i) separates points of $\Lam$, then the evaluation map $\ev:\Lam\rightarrow{(D^0)^\star=(\{\partial_\theta\mid\theta\in\uu\})^\star}$ is injective. If moreover the image of $\ev$ is discrete, then $D_{(\u,\u',\Gamma')}$ is permissible. 
	
	(iii) {Here we give an example of a subspace $\u$ of $\mathrm{Hom}_\bbbz(\Lam,\bbbk)$ which separates points of $\Lam$.} 
	Assume that $\Lam=\sum_{i=1}^n\bbbz\sg_i$, then $\mathrm{Hom}_\bbbz(\Lam,\bbbk)=\sum_{i=1}^n\bbbk\theta_i$, where
	$\theta_i(\sg_j)=\d_{i,j}$.  If 
	$p_1,\ldots, p_n$ are any $\bbbz$-linearly independent elements of $\bbbk$, then
	the one dimensional space $\bbbk(p_1\theta_1+\cdots+p_n\theta_n)$ separates points of $\Lam$.}
	
	(iv) {Let $\u=\mathrm{Hom}_\bbbz(\Lam,\bbbk)$.  Then $\text{ev}:\Lam\rightarrow \dd^\star$ is a discrete embedding. Thus $D_{(\u,\u^\prime,\Gamma^\prime)}$ is a  permissible {subalgebra} of $\mathrm{SCDer}_\bbbk(\ll)$. In particular, $D_{(\u,\u,\Gamma)}=\mathrm{SCDer}_\bbbk(\ll)$ and $D_{(\u,\{0\},\{0\})}=\dd$ are permissible subalgebras of $\mathrm{SCDer}_\bbbk(\ll)$.}
\end{exa}

\begin{parag}\label{parag3}
Let $\ll$ be a centerless Lie torus, $D$ be a permissible subalgebra of
$\scd(\ll)$ and $\kappa$ be an affine cocycle on $D$, see  \pref{parag2}. Let
$$E=E(\ll,D,\kappa):=\ll\oplus D^{gr*}\oplus D$$
be the Lie algebra with the Lie bracket
\begin{eqnarray*}
	[x_1+c_1+d_1,x_2+c_2+d_2]&=&([x_1,x_2]_\ll+d_1(x_2)-d_2(x_1))\\
	&+&(\sg_D(x_1,x_2)+d_1.c_2-d_2.c_1+\kappa(d_1,d_2))\\
	&+&[d_1,d_2]
\end{eqnarray*}
for $x_1,x_2\in \ll,c_1,c_2\in D^{gr*},d_1,d_2\in D$, where $[\:,\:]_\ll$ denotes the Lie bracket of $\ll$, $[d_1,d_2]=d_1d_2-d_2d_1$, and $\sg_D:\ll\times\ll\rightarrow D^{gr*}$ is defined by
$$\sg_D(x,y)(d)=(d(x)|y)\;\text{for all}\;x,y\in \ll,d\in D.$$
In fact $\sg_D$ is a $2$-cocycle for $\ll$ whose values are in the trivial $\ll$-module $D^{gr*}$. Also $\sg_D$ respects the gradings of $\ll$ and $D^{gr*}$. It follows that $E$ has a non-degenerate invariant symmetric bilinear form $\fm$ given by
\begin{equation}\label{ira1}
(x_1+c_1+d_1,x_2+c_2+d_2)=(x_1,x_2)_\ll+c_1(d_2)+c_2(d_1).
\end{equation}
\end{parag}
Let 
\begin{equation}\label{tt}
	H:=\fh\oplus(D^0)^\star\oplus D^0,\hbox{ where }\fh:=\ll_0^0,
\end{equation}
see Definition \ref{def1}. Then we have the following.

\begin{thm}\cite[{Theorem 16}]{Neh04}
	(a) The algebra $E(\ll,D,\kappa)$ constructed in \pref{parag3} is an extended affine Lie algebra of nullity $n$ with respect to the form (\ref{ira1}) and the Cartan subalgebra $H$ given in (\ref{tt}). Moreover, $E_c=\ll\oplus  {D^{gr}}^\star$.
		
	(b) Conversely, let $E$ be an extended affine Lie algebra, and let $\ll=E_{cc}$ be its centerless core of type $(\Lam,\Delta)$. Then there exists a unique permissible subalgebra $D\sub\scd(\ll)$ and an affine cocycle $\kappa$ on $D$ such that $D$ induces the $\Lam$-grading of $\ll$, and $E\cong E(\ll,D,\kappa)$. 
\end{thm}


\section{\bf Chevalley and Pre-Chevalley involutions}\setcounter{equation}{0}
We discuss the possibility of lifting a Chevalley involution on a centerless Lie torus $\ll$ to the extended affine Lie algebra $E(\ll,D,\kappa),$ where
$D$ is a permissible subalgebra of skew centroidal derivations of $\ll$ and $\kappa$ is an affine cocycle.
Let $\ll$ be a Lie $\Lam$-torus of type $\Delta$ as in Definition \ref{def1}. 


\subsection{Induced pre-Chevalley involutions}\label{Induces pre-Chevalley involutions} We begin with the formal definition of a Chevalley involution for a centerless Lie torus.

\begin{DEF}\label{1def}
	{{(i) We call an involution (an automorphism of order $2$) $\tau$ of a $\Lam$-graded algebra $\aa$, 
	a {\it pre-Chevalley involution} if $\tau(\aa^\lam)=\aa^{-\lam}$ for all $\lam\in\Lam.$ }}
	
	{(ii) A pre-Chevalley involution $\tau$ of $\ll$ is called a {\it Chevalley involution} if  $\tau(h)=-h$ for $h\in\ll^0_0$.}
	\end{DEF}

\begin{lem}
	Assume that $\tau$ is a Chevalley involution for $\ll$. Then $\tau(\ll_\a)=\ll_{-\a}$ for $\a\in\Delta$.
\end{lem}
\proof
First consider the root-grading pair $(\fg,\fh)$ for $\ll$ where $\fg$ is the subalgebra of $\ll$ generated by $\{\ll_\a^0\}_{\a\in\Delta^\times}$ and $\fh=\sum_{\a\in\Delta^\times}[\ll_\a^0,\ll_{-\a}^0]$. Then
$$\ll_\a=\{x\in\ll\mid[h,x]={\a}(h)x\;\text{for}\;h\in\fh\},$$
for $\a\in Q$, {where $\a$ is identified in the natural way as an element of $\fh^\star$}, see \cite[Proposition 1.2.2]{ABFP09}. Note that $\fh\subseteq\ll_0^0$ and since $\tau_{|_{\ll_0^0}}=-\id$ we have
$$[h,\tau(x)]=\tau([-h,x])=-{\a}(h)\tau(x),$$
for $\a\in Q,x\in\ll_\a$ and $h\in\fh$. So $\tau(\ll_\a)=\ll_{-\a}$, for $\a\in Q$.\qed

\begin{rem}\label{rem4}
	{Assume that $\tau$ is a pre-Chevalley involution for $\ll$. It is easy to see that $\tau\partial_\theta=-\partial_\theta\tau$
		for $\theta\in\mathrm{Hom}_\bbbz(\Lam,\bbbk)$.}
\end{rem}

Suppose that $\tau$ is a pre-Chevalley involution for $\ll$.
We consider associative algebra automorphisms (both denoted by $\bar\tau$ for the sake of notation),
\begin{equation}\label{auto1}
	\begin{array}{c}
		\bar\tau:\End(\ll)\rightarrow\End(\ll)\\\phi\mapsto\tau\phi\tau^{-1}
	\end{array}\andd
	\begin{array}{c}
		\bar\tau:\End(\ll)^\star\rightarrow\End(\ll)^\star\\
		\chi\mapsto\bar\tau(\chi):\phi\mapsto{\chi(\tau^{-1}\phi\tau)}
	\end{array}
\end{equation}
for $\phi\in\End(\ll)$ and $\chi\in\End(\ll)^\star$. {Then $\bar\tau$ can also be considered as an involutions on $\gl(\ll)$, or $\gl(\ll)^\star$.}

{The $\Lam$-grading of $\ll$ induces a $\Lam$-gradings on
$$\mathrm{gr}\End(\ll):=\sum_{\lam\in\Lam}\End(\ll)^\lam$$ where
$$\End(\ll)^\lam=\{\phi\in\End(\ll)\mid\phi(\ll^\mu)\subseteq\ll^{\lam+\mu},\;\text{for all}\;\mu\in\Lam\}.$$
We also consider the graded dual $\mathrm{gr}{\End(\ll)^{gr}}^\star=\sum_{\lam\in\Lam}(\End(\ll)^\lam)^\star$ of $\mathrm{gr}\End(\ll)$ as a $\Lam$-graded algebra with 
$$(\mathrm{gr}{\End(\ll)^{gr}}^\star)^\lam=(\End(\ll)^{-\lam})^\star,$$
for $\lam\in\Lam$. Under the commutator bracket, the resulting $\Lam$-graded Lie algebras will be denoted by
$\mathrm{gr}\gl(\ll)$ and $\mathrm{gr}{\gl(\ll)^{gr}}^\star$, respectively. 

{Next, for $\mu\in\Gamma$ and $\lam\in\Lam$, we define
\begin{equation}\label{tt1}
c^{(\mu)}_\lam(\chi^{-\mu}\partial_\theta):=\theta(\lam),
	\end{equation}
	and extend it to an element of ${\scd(\ll)^{gr}}^\star$ by $c^{(\mu)}_\lam(\chi^\nu\partial_\gamma)=0$ for $\nu\not=-\mu.$
Using these facts, now we have the following result.}

\begin{lem}\label{sky23} Let $\tau$ be a pre-Chevalley involution on $\ll$, and consider the automorphism $\bar\tau$ as in (\ref{auto1}).
	
	{(i) $\bar\tau$ restricts to pre-Chevalley involutions on $\mathrm{gr}\End(\ll)$ and $\mathrm{gr}{\End(\ll)^{gr}}^\star$.
		
		(ii) $\bar\tau$ restricts to a pre-Chevalley involution on $\cent_\bbbk(\ll)$. In particular, $\cent_\bbbk(\ll)$ has a basis $\{\chi^\mu\mid\mu\in\Gamma\}$ shch that $\bar\tau(\chi^\mu)=\chi^{-\mu}$, $\mu\in\Gamma$.}
	\end{lem}

\proof
{(i) For $\phi\in\End(\ll)^\lam$, we have
	$$\bar\tau(\phi)(\ll^\mu)=\tau\phi\tau^{-1}(\ll^\mu)=\tau\phi(\ll^{-\mu})\subseteq\tau(\ll^{\lam-\mu})= \ll^{-\lam+\mu},$$
	for all $\mu\in\Lam$. So $\bar\tau(\phi)\in\End(\ll)^{-\lam}$. This shows that $\bar\tau$ is a pre-Chevalley involution for $\mathrm{gr}\End(\ll)$. Similarly {$\bar\tau({(\End(\ll)^\lam)}^\star)={(\End(\ll)^{-\lam})}^\star$}, for all $\lam\in\Lam$. Thus $\bar\tau$ is a pre-Chevalley involution for $\mathrm{gr}\End(\ll)^\star$.} 

{(ii) Let $\chi\in\cent_\bbbk(\ll)$, then we have
	$$\tau\chi\tau^{-1}([x,y])=[x,\tau\chi\tau^{-1}(y)]=[\tau\chi\tau^{-1}(x),y],\;\forall x,y\in \ll.$$
	So $\bar\tau(\chi)=\tau\chi\tau^{-1}\in\cent_\bbbk(\ll)$. Also note that $\cent_\bbbk(\ll)^\lam=\End(\ll)^\lam\cap\cent_\bbbk(\ll)$. Thus $\bar\tau$ restricts to a a pre-Chevalley involution for $\cent_\bbbk(\ll)$.}

{Now recall that
	$\cent_\bbbk(\ll)=\bigoplus_{\mu\in\Gamma}\bbbk\chi^\mu,$
	where $\chi^\mu$ acts on $\ll$ as an endomorphism of degree $\mu$, and $\Gamma$ is the central grading group of $\ll$. So $\bar\tau(\chi^\mu)=\eta_\mu\chi^{-\mu}$, for some $\eta_\mu\in\bbbk^\times$.
Since $\bar\tau$ has order $2$, $\eta_\mu\eta_{-\mu}=1.$ Then
for each $\mu$,
$$\bar\tau(\eta_\mu^{\frac{-1}{2}}\chi^\mu)
=\eta_\mu\eta_\mu^{\frac{-1}{2}}\chi^{-\mu}
={\eta_\mu^\frac{1}{2}}\chi^{-\mu}=
\eta_{-\mu}^{\frac{-1}{2}}\chi^{-\mu}.$$
Thus replacing each $\chi^\mu$ by $\eta_\mu^{\frac{-1}{2}}\chi^\mu$, we may assume that $\bar\tau(\chi^\mu)=\chi^{-\mu}$.\qed
}
	
	Using Lemma \ref{sky23}, {\it we assume from now on} that 
\begin{equation}\label{sky24}\cent_\bbbk(\ll)=\sum_{\mu\in\Gamma}\bbbk\chi^\mu\hbox{ with }\bar\tau(\chi^\mu)=\chi^{-\mu}.
\end{equation}

\begin{lem}\label{llem1} Let $\tau$ be a pre-Chevalley involution on $\ll$, and $\bar\tau$ be as in (\ref{auto1}).

	{(i) $\bar\tau$ induces a pre-Chevalley involution on $\mathrm{gr}\gl(\ll)$. Moreover $\bar\tau(\chi^\mu\partial_\theta)=\chi^{-\mu}\partial_{-\theta}$ for each $\mu$ and $\theta$. In particular, $\bar\tau(\chi^\mu\dd)=\chi^{-\mu}\dd$, $\mu\in\Gamma$, and $\bar\tau$ restricts to a pre-Chevalley involution on $\mathrm{SCDer}(\ll)$.}

{(ii)
$\bar\tau$ induces  a pre-Chevalley involution on  $\mathrm{gr}{\gl(\ll)^{gr}}^\star$. Moreover $\bar\tau(c^{(\mu)}_\lam)=c^{(-\mu)}_{-\lam}$ for each $\mu$ and $\lam$. In particular,  $\bar\tau$ restricts to a pre-Chevalley involution on
{${\scd(\ll)^{gr}}^\star$}.}

\end{lem}
\proof
	{(i) For $\theta\in\mathrm{Hom}_\bbbz(\Lam,\bbbk)$ and $x^\lam\in\ll^\lam$, we have
		$$\bar\tau(\partial_\theta)(x^\lam)=\tau\partial_\theta\tau(x^\lam)=-\theta(\lam)\tau^2(x^\lam)=\partial_{-\theta}(x^\lam).$$
	Thus $\bar\tau(\partial_\theta)=\partial_{-\theta}.$
	Now the statement follows from this and Lemma \ref{sky23}.}

{(ii)
	It follows that $\{c^{(\mu)}_\lam\mid\mu\in\Gamma,\;\lam\in\Lam\}$ spans
	${\scd(\ll)^{gr}}^\star$, see \cite[Proposition 5.2.4]{Nao10}. 
	Now for any $\theta$ and $\nu$,
	\begin{eqnarray*}
		\bar\tau(c_\lam^{(\mu)})(\chi^\nu\partial_\theta)=
		c_\lam^{(\mu)}(\bar\tau(\chi^\nu\partial_\theta))
		=c_\lam^{(\mu)}(\chi^{-\nu}\partial_{-\theta})
		=
	 c_{-\lam}^{(-\mu)}(\chi^\nu\partial_\theta).
	\end{eqnarray*}
	Thus
	$\bar\tau(c^{(\mu)}_\lam)=c^{(-\mu)}_{-\lam}\in ({\scd(\ll)^{gr}}^\star)^{-\mu}.$
{That is,}
	$\bar\tau$ maps {$({\scd(\ll)^{gr}}^\star)^{\mu}$} onto {$({\scd(\ll)^{gr}}^\star)^{-\mu}$}, and $\bar\tau$ restricts to a pre-Chevalley involution for {${\scd(\ll)^{gr}}^\star$}, as a subalgebra of $\mathrm{gr}{\gl(\ll)^{gr}}^\star$.}}\qed

{\begin{rem}\label{sky16} Let $D=\sum_{\mu\in\Gamma}D^\mu$ be a permissible subalgebra of $\scd(\ll)$ and $C={D^{gr}}^\star$.
	Consider the maps $c^{(\mu)}_\lam$, $\mu\in\Gamma$, given in (\ref{tt1}).
	For simplicity of notation, we denote
	the restriction of $c_\lam^{(\mu)}$ to $D^{-\mu}$ by $c^{(\mu)}_\lam$, again. Then it follows that $C^\mu=\span_{\bbbk}\{c^{(\mu)}_\lam\mid\lam\in\Lam\}$. 
\end{rem}
}

  \subsection{Induced affine cocycles}\label{Induced $2$-cocycles}
 { Let $\ll$ be a centerless Lie torus, 
 	$D=\bigoplus_{\mu\in\Gamma}D^\mu$ be a {permissible} subalgebra of  $\mathrm{SCDer}_\bbbk(\ll)$, and  $\kappa:D\times D\rightarrow D^{gr*}$ be an affine cocycle. Assume that $\tau$ is a pre-Chevalley involution for the Lie torus  $\ll$.}
 {We set
 \begin{equation}\label{sky14}
 D_\tau:=\bar\tau(D)\andd D^\mu_\tau:=\bar\tau(D^\mu),\;\mu\in\Gamma.
 \end{equation}
Then we get the natural identification $\bar\tau(D^{gr^\star})=D_\tau^{gr^\star}=\sum_{\mu\in\Gamma}{D_\tau^\mu}^\star.$ We next set
\begin{equation}\label{sky15}
 \begin{array}{c}
 \kappa_\tau:D_\tau\times D_\tau\rightarrow D_\tau^{gr^\star}\\
 \kappa_{\tau}(\bar\tau(d),\bar\tau(d'))=\bar\tau\kappa(d,d').
 \end{array}
\end{equation}
 }

 \begin{lem}\label{1lem} {(i) $D_\tau$ is  a permissible subalgebra of $\mathrm{SCDer}(\ll)$.}
 	
 	{(ii) 
{The  map $\kappa_\tau$
 defines an affine cocycle on $D_\tau$.}}
 \end{lem}
 
 \proof { (i) is clear as  $\bar\tau(\partial_\theta)=-\partial_\theta$ and so
 $D_\tau^0=D^0.$  }
 
(ii) We first check the $\circlearrowleft$ part in the definition of an affine cocycle. 
For $d_1,d_2,d_3\in D$,
		\begin{eqnarray*}
			\sum_{(i,j,k)\circlearrowleft}\kappa_\tau([\bar\tau(d_i),\bar\tau(d_j)],\bar\tau(d_k))&=&\sum_{(i,j,k)\circlearrowleft}\bar\tau(\kappa([d_i,d_j],d_k)\\
			&=&\bar\tau(\sum_{(i,j,k)\circlearrowleft}\kappa([d_i,d_j],d_k))\\
			&=&\bar\tau(\sum_{(i,j,k)\circlearrowleft}d_i\cdot\kappa (d_j,d_k))\\
			&=&\sum_{(i,j,k)\circlearrowleft}\bar\tau(d_i\cdot\kappa (d_j,d_k)).
		\end{eqnarray*}
	Now for $d\in D$, we have
		\begin{eqnarray*}
			\sum_{(i,j,k)\circlearrowleft}	\bar\tau(d_i\cdot\kappa (d_j,d_k))(d)&=&\sum_{(i,j,k)\circlearrowleft}d_i\cdot\kappa(d_j,d_k)(\tau^{-1}d\tau)\\
			&=&\sum_{(i,j,k)\circlearrowleft}\kappa(d_j,d_k)([\tau^{-1}d\tau,d_i])\\
			&=&\sum_{(i,j,k)\circlearrowleft}\kappa(d_j,d_k)(\tau^{-1}[d,\tau d_i\tau^{-1}]\tau)\\
			&=&\sum_{(i,j,k)\circlearrowleft}\bar\tau(\kappa(d_j,d_k))([d,\bar\tau(d_i)])\\
			&=&\sum_{(i,j,k)\circlearrowleft}\bar\tau(d_i)\cdot\bar\tau(\kappa(d_j,d_k))(d).
		\end{eqnarray*}
		So $\sum_{(i,j,k)\circlearrowleft}	\bar\tau(d_i\cdot\kappa (d_j,d_k))=\sum_{(i,j,k)\circlearrowleft}\bar\tau(d_i)\cdot\bar\tau(\kappa(d_j,d_k))$. This together with the last equality gives
		\begin{eqnarray*}
			\sum_{(i,j,k)\circlearrowleft}\kappa_\tau([\bar\tau(d_i),\bar\tau(d_j)],\bar\tau(d_k))&=&\sum_{(i,j,k)\circlearrowleft}\bar\tau(d_i)\cdot\bar\tau(\kappa(d_j,d_k))\\
			&=&\sum_{(i,j,k)\circlearrowleft}\bar\tau(d_i)\cdot\kappa_\tau(\bar\tau(d_j),\bar\tau(d_k)),
		\end{eqnarray*}
so the proof of $\circlearrowleft$ is completed.
		
		Next, for $\mu_1,\mu_2\in\Gamma$, since $\kappa(D^{-\mu_1},D^{-\mu_2})\subseteq(D^{\mu_1+\mu_2})^\star$ we have
		\begin{eqnarray*}		\kappa_\tau(D_\tau^{\mu_1},D_\tau^{\mu_2})&=&\kappa_\tau(\bar\tau(D^{-\mu_1}),\bar\tau(D^{-\mu_2}))\\
			&=&\bar\tau(\kappa(D^{-\mu_1},D^{-\mu_2}))\\
			&\subseteq&\bar\tau((D^{\mu_1+\mu_2})^\star)=(D_\tau^{-\mu_1-\mu_2})^\star.
		\end{eqnarray*}
		To check that the remaining properties of an affine cocycle hold for $\kappa_\tau$ is straightforward.
\qed

We need the following concept in the sequel.
\begin{DEF}\label{def2}
	{Let $\ll$ be a Lie torus, $D$ be  a permissible subalgebra of $\mathrm{SCDer}(\ll)$ and $\kappa$ be an affine cocycle on $D$. We call $D$, $\tau$-{\it invariant} if $D_\tau=D$ and call the  pair $(D,\kappa)$,  {\it $\tau$-invariant} if $(D_\tau,\kappa_\tau)=(D,\kappa)$.}
\end{DEF}	
\begin{rem}\label{sky33}
{Suppose $D$, $\kappa$ and $\tau$ are as in Definition \ref{def2}, and suppose that
$D=D_\tau$. We know that $\kappa=\kappa_\tau$ if and only if for any
$d_1,d_2,d_3\in D$,
$$\begin{array}{c}
	\bar\tau\kappa(d_1,d_2)(d_3)=\kappa(\bar\tau(d_1),\bar\tau(d_2))(d_3)\\
	\Longleftrightarrow\\
	\kappa(d_1,d_2)(\bar\tau(d_3))=\kappa(\bar\tau(d_1),\bar\tau(d_2))(d_3).
	\end{array}
	$$
Then since $\kappa(D^\mu,D^\nu)\sub ({D^{-\mu-\nu}})^\star$, for all $\mu$, $\nu$, we get $\kappa=\kappa_\tau$ if and only if
$$\kappa(\chi^{\mu}\partial_\theta,\chi^{\nu}\partial_\gamma)(\chi^{-\mu-\nu}\partial_{-\eta})
=\kappa(\chi^{-\mu}\partial_\theta,\chi^{-\nu}\partial_\gamma)(\chi^{\mu+\nu}\partial_\eta),
$$
for any $\mu,\nu,\theta,\gamma$ and $\eta$ with $\theta(\mu)=\gamma(\nu)=\eta(\mu+\nu)=0$.}
\end{rem}

Let $D=\sum_{\mu\in\Gamma}D^\mu$ be a permissible subalgebra of $\scd(\ll)$. We set
\begin{equation}\label{sky20}
\u_D^\mu:=\{\theta\in\mathrm{Hom}_\bbbz(\Lam,\bbbk)\mid\chi^\mu\partial_\theta\in D^\mu\}.
\end{equation}

\begin{lem}\label{sky19}
A permissible subalgebra $D$ of $\scd(\ll)$ is $\tau$-invariant if and only if $\u_D^\mu=\u_D^{-\mu}$ for all $\mu\in\Gamma$. {In particular, if $D=D^0$ or $D=\scd(\ll)$, then $D$ is $\tau$-invariant.}
\end{lem}

\proof Using Lemma \ref{llem1}(i), we have 
\begin{eqnarray*}
\bar\tau(D)=D&\Longleftrightarrow&\bar\tau(D^\mu)=
D^{-\mu}\hbox{ for all }\mu\in\Gamma,\\
&\Longleftrightarrow&
\bar\tau(\chi^\mu\partial_\theta)=\chi^{-\mu}\partial_{-\theta}\in D^{-\mu}\hbox{ for all }
\mu\in\Gamma,\;\theta\in\u_D^\mu,\\
&\Longleftrightarrow&
\u_D^\mu=\u_D^{-\mu}\hbox{ for all }\mu\in\Gamma.
\end{eqnarray*}
\qed

In part (ii) of the below example we construct a class of permissible subalgebras of $\scd(\ll)$ whose elements are not
$\tau$-invariant, with respect to any $\tau$.

\begin{exa}\label{sky-10} Assume $\ll$ is a centerless Lie torus with centroidal grading group $\Gamma\sub\Lam$.
	
		(i) Consider the pair $(D,\kappa)$, where $D$ is permissible and  $\kappa$ is an affine cocycle. Then $D\cap D_\tau=D^0+\sum_{\mu\in\Gamma\setminus\{0\}}(D^\mu\cap D^{-\mu}_\tau)$ is permissible, $\kappa+\kappa_\tau$ is an affine cocycle on $D\cap D_\tau$ and  $(D\cap D_\tau,\kappa+\kappa_\tau)$ is a $\tau$-invariant pair.
	
	{ (ii) Let rank $\Lam>1$. We fix a triple $(\u,\u_+ , \u_{-})$ of subspaces of $\hbox{Home}_\bbbz(\Lam,\bbbk)$. Also, we fix $0\not=\gamma\in\Gamma$ {and for $\mu\in\Gamma$ set}
		$${\u^\mu}=\left\{\begin{array}{ll}
			\u&\hbox{if }\mu=0,\\
			\{\theta\in\u_\pm\mid \theta(\mu)=0\}&\hbox{if }\mu=\pm\gamma,\\
			\{0\}&\hbox{otherwise},
		\end{array}
		\right.
		$$
		and
		$$D^\mu:=\chi^\mu\{\partial_\theta\mid\theta\in{\u^\mu}\}.
		$$			
		Then $D:=\sum_{\mu\in\Gamma}D^\mu$ is a $\Gamma$-graded subspace of $\scd(\ll)$.
		Since by (\ref{eq7}), $[D^{\pm\gamma},D^{\pm\gamma}]=\{0\}$, it follows that  $D$ is a $\Gamma$-graded subalgebra. 
		Note that if we take $\u=\hbox{Hom}_\bbbz(\Lam,\bbbk)$ then
		$D$ is a permissible subalgebra of $\scd(\ll)$
		{with $\u^\mu_D=\u^\mu$ for each $\mu\in\Gamma$, see (\ref{sky20})}. Further,
		we choose $\u_+,\u_-$ such that {$\u^\gamma\not=\u^{-\gamma}$, 
			for instance}  $\u_+=\hbox{Hom}_\bbbz(\Lam,\bbbk)$ and $\u_{-}=\{0\}$.  Then by Lemma \ref{sky19}, $D$ is not $\tau$-invariant, with respect to any Chevalley involution $\tau$.
	}
\end{exa}

\subsection{Induced Chevalley involutions on EALAs}
{ We proceed with the same notation and assumptions as in the preceding sections. In particular, $\tau$ is a pre-Chevalley involution for the Lie torus $\ll$ and $\bar\tau$ is as in (\ref{auto1}). In this subsection, we {discuss} the natural extension of $\bar\tau$ to $\ll\oplus\End(\ll)\oplus\End(\ll)^\star$, denoted again by $\bar\tau$, see Theorem \ref{induce} below.
}

\begin{DEF}\label{sky-11} {Assume that $(E,\fm,\hh)$ is an extended affine Lie algebra with root system $R$. We call an involution $\tau$ of $E$ a {\it Chevalley involution} if $\tau(h)=-h$ for all $h\in\hh$.}
\end{DEF}


\begin{thm}\label{induce}
Let $\ll$ be a centerless Lie $\Lam$-torus of type $\Delta$, $D$ be a permissible subalgebra of $\scd(\ll)$ and $\kappa$ be an affine cocycle on $D$. Let  $\tau$  be a pre-Chevalley involution for $\ll$, {and $\bar\tau$ be the induced vector space isomorphism given by}
$$\begin{array}{c}
\bar\tau:\ll\oplus\End(\ll)\oplus\End(\ll)^\star{\longrightarrow \ll\oplus\End(\ll)\oplus\End(\ll)^\star}\\
x+\phi+\Psi\mapsto\tau(x)+\bar\tau(\phi)+\bar\tau(\Psi),
\end{array}
$$
see (\ref{auto1}). Let $E=E(\ll,D,\kappa)$ {and $E^\tau:=E(\ll,D_\tau,\kappa_\tau)$}. Then

(i) {$\bar\tau$ restricts to a Lie algebra isomorphism from $E_c=\ll\oplus D^{gr^\star}$ onto $E^\tau_c=\ll\oplus D_\tau^{gr^\star}$. In particular if $D$ is $\tau$-invariant, then $\bar\tau$ restricts to a pre-Chevalley involution for $E_c$.}

{(ii) {$\bar\tau$ restricts to a Lie algebra isomorphism from {$E$ onto $E^\tau$}.} In particular if {$(D,\kappa)$} is $\tau$-invariant, then $\bar\tau$ restricts to a pre-Chevalley involution for $E$.}

If in addition $\tau$ is a Chevalley involution, then so is $\bar\tau.$
\end{thm}

\proof
We first note that by Lemma \ref{1lem}, $E^\tau$ is an extended affine Lie algebra. We only need to show that $\bar\tau$ is a Lie algebra homomorphism.
Let $x_1,x_2\in\ll$, $c_1,c_2\in C$ and $d_1,d_2\in D$.
Then
\begin{equation*}
	\begin{split}
		\bar\tau[x_1+c_1+d_1,x_2+c_2+d_2]=&[\tau(x_1),\tau(x_2)]_\ll+\tau(d_1(x_2))-\tau(d_2(x_1))\\
		+&\bar\tau(\sg_D(x_1,x_2)+\bar\tau (d_1.c_2)-\bar\tau (d_2.c_1)+\bar\tau\kappa(d_1,d_2))\\
		+&[\bar\tau(d_1),\bar\tau(d_2)].
	\end{split}
\end{equation*}
Now $\tau(d_1(x_2))=\tau d_1\tau^{-1}(\tau(x_2))=\bar\tau(d_1)(\tau(x_2))$. Similarly, we have $\tau(d_2(x_1))=\tau d_2\tau^{-1}(\tau(x_1))=\bar\tau(d_2)(\tau(x_1))$. Also for $d\in D$,
\begin{eqnarray*}
\bar\tau\sg_D(x_1,x_2)(d)&=&\sg_D(x_1,x_2)(\tau^{-1}d\tau)\\
&=&
(\tau^{-1}d\tau(x_1),x_2)\\
&=&(d(\tau(x_1),\tau(x_2))\\
&=&\sg_D(\tau(x_1),\tau(x_2))(d).
\end{eqnarray*}
Therefore $\bar\tau\sg_D(x_1,x_2)=\sg_D(\tau(x_1),\tau(x_2)).$ Next for $d\in D$, 
\begin{eqnarray*}
\bar\tau(d_1\cdot c_2)(d)&=&d_1\cdot c_2(\tau^{-1}d\tau)
\\&=&c_2[\tau^{-1}d\tau,d_1]\\
&=&
c_2\tau^{-1}[d,\tau d_1\tau^{-1}]\tau\\
&=&\bar\tau(d_1)\cdot\bar\tau(c_2)(d).
\end{eqnarray*}
Thus 
$\bar\tau(d_1\cdot c_2)=\bar\tau(d_1)\cdot\bar\tau(c_2).$ 
Also $\bar\tau\kappa(d_1,d_2)=\kappa_\tau(\bar\tau(d_1),\bar\tau(d_2))$ for $d_1,d_2\in D$. So $\bar\tau$ is a Lie algebra isomorphism from $E(\ll,D,\kappa)$ onto $E(\ll,D_\tau,\kappa_\tau)$. {In particular $\bar\tau$ restricts to a Lie algebra isomorphism from $E_c$ onto $E^\tau_c$.}

{Next, if $(D,\kappa)$ is $\tau$-invariant, i.e., $D=D_\tau$ and $\kappa=\kappa_\tau$, then clearly $\bar\tau$ is a pre-Chevalley involution for $E$. {In particular, if $D$ is $\tau$-invariant, i.e., $D=D_\tau$, then $\bar\tau$ is a pre-Chevalley involution for $E_c$.} The last part of the statement is now clear.}\qed

\begin{cor}\label{sky31}
	{Let $E=E(\ll,D,\kappa)$ where $\ll$ is a centerless Lie torus, and $D=D^0$ or $D=\scd(\ll)$. Then any  Chevalley involution for $\ll$ extends to  a Chevalley involution for $E_c$.}
\end{cor}
\proof {Since by Lemma \ref{sky19} both $D^0$ and $\scd(\ll)$ are $\tau$-invariant, with respect to any Chevalley involution $\tau$ for $\ll$, the statement follows from Theorem \ref{induce}.\qed}
 
\begin{rem}\label{remneg}
{(i) In order to study the concept of ``modular theory" for extended affine Lie algebras, in \cite{AFI22} the core of a reduced extended affine Lie algebra $E$ of rank $>1$ is equipped with an integral structure.  This in fact is achieved by assuming that $E_c$ is equipped with a Chevalley involution. Now suppose $E=E(\ll,D,\kappa)$ and suppose that the centerless Lie torus $\ll$ is equipped with a Chevalley involution $\tau$. If $D$ is $\tau$-invariant, then Theorem \ref{induce}(i) guarantees the existence of a Chevalley involution for the core. Special cases $D=D^0$ and $D=\scd(\ll)$ which are $\tau$-invariant with respect to any Chevalley involution $\tau$ for $\ll$ are discussed in Corollary \ref{sky31}. Much of the most interesting examples of extended affine Lie algebras {appeared} in the literature satisfy conditions of Corollary \ref{sky31}. In fact, as far as the realization of extended affine Lie algebras is concerned, in the literature we are not aware of any specific example of an extended affine Lie algebra that does not fall in the conditions of Corollary \ref{sky31}, except for the one we gave in Example \ref{sky-10}. The existence of Chevalley involutions for the centerless Lie tori will be discussed in Section \ref{Chevalley involutions for centerless Lie tori}.}

{(ii) Concerning Theorem \ref{induce}, it is a demanding question to ask under which conditions an affine cocycle $\kappa$ is $\tau$-invariant, namely $\kappa_\tau=\kappa$. In  this regard, it is natural to first investigate the situation for the known affine cocycles.
Suppose  $\ll=\fg\otimes\aa$, where $\fg$ is a finite dimensional simple Lie algebra
of simply laced type and $\aa$ is the algebra of Laurent polynomials in $\nu$ variables, and consider $D=\scd(\ll)$.  
In \cite[Remark 3.71]{BGK96} a non-trivial affine cocycle
$\kappa:D\times D\rightarrow {D^{gr}}^\star$ is provided (see also \cite{EM94}). In \cite[Appendix II, Proposition 8]{Kry00}, it is shown that $\kappa$ is $\tau$-invariant for any {$\tau\in\Aut(\cent(\ll))$}. Thus $(D,\kappa)$ is $\tau$-invariant and so Theorem \ref{induce}(ii) applies.  
}

{(iii) {In Theorem \ref{induce}, assume that $D$ is $\tau$-invariant. Then by  restricting the involution $\bar\tau$  to $D\oplus {D^{gr}}^\star$, we conclude that the two short exact sequences
$$
\begin{array}{c}
0\rightarrow {D^{gr}}^\star\rightarrow D\oplus_\kappa {D^{gr}}^\star\rightarrow D\rightarrow 0,\vspace{2mm}\\
0\rightarrow {D^{gr}}^\star\rightarrow D\oplus_{\kappa_\tau} {D^{gr}}^\star\rightarrow D\rightarrow 0,
\end{array}
$$
are isomorphic as extensions of $D$ by the $D$-module ${D^{gr}}^\star,$ meaning that $\kappa$ and $\kappa_\tau$ represent the same class in $H^2(D,{D^{gr}}^\star)$ (see \cite[Theorem 7.6.3]{Wei94}).}}
\end{rem}

\section{\bf Pre-Chevalley involution for $\Lam$-Tori}\setcounter{equation}{0}
\label{tori}

The notion of a $\Lam$-torus arises in the study of Lie tori. Roughly speaking, a centerless Lie $\Lam$-torus can be constructed as a matrix Lie algebra coordinatized by a $\Lam$-torus. Depending on the type, it is a unital associative, Jordan or alternative {algebra.} We begin in this section by recalling the definition of a $\Lam$-torus. Then we review the classification results of $\Lam$-tori and using theses results we show that any $\Lam$-torus admits a pre-Chevalley involution. Throughout the section, $\Lam$ is a free abelian group of rank $n$ and by an algebra over $\bbbk$ we mean a unital Jordan algebra, alternative algebra or associative algebra over $\bbbk$. Assuming $\aa$ and $\bb$ are $\Lam$-graded algebras, then we say that $\aa$ and $\bb$ are graded-isomorphic, if there is an algebra isomorphism $\varphi:\aa\rightarrow\bb$ that preserves the $\Lam$-grading.

\begin{DEF}\label{def7}
	Let $\aa=\bigoplus_{\lam\in\Lam}\aa^\lam$ be a $\Lam$-graded algebra over $\bbbk$. $\aa$ is called a {\it $\Lam$-torus} or {\it $n$-torus} if:
	
	(i) all $0\neq x\in \aa^\lam$ are invertible and $\dim_\bbbk \aa^\lam\leq1$, for all $\lam\in\Lam$,
	
	(ii) $\supp_\Lam(\aa)$ generates $\Lam$, where $\supp_\Lam(\aa)=\{\lam\in\Lam\mid \aa^\lam\neq\{0\}\}$.\\
	If a $\Lam$-torus $\aa$ is a Jordan algebra, alternative algebra or associative algebra then $\aa$ is called a {\it Jordan, alternative, or associative $\Lam$-torus ($n$-torus)}, respectively.
\end{DEF}

{To have a pre-Chevalley involution $\tau$ of a $\Lam$-torus $\aa=\sum_{\a\in\supp(\aa)}\aa^\lam$, it looks natural to first consider an appropriate basis $\{x^\lam\}_{\lam\in\supp(\aa)}$ for $\aa$, $x^\lam\in \aa^\lam$, and then define $\tau(x^\lam)=x^{-\lam}$, $\lam\in\supp(\aa)$. This gives a vector space isomorphism which maps $\aa^\lam$ onto $\aa^{-\lam}$. Therefore, we get a pre-Chevalley involution for the $\Lam$-torus $\aa$, if we show that $k(\a,\b)=k(-\a,-\b)$ for $\a,\b\in\Lam$, where
	$x^\a x^\b=k(\a,\b)x^{\a+\b}$, $k(\a,\b)\in \bbbk^\times.$} 

\begin{thm}\label{thm-2}
	{Let $\aa$ be a Jordan, alternative, or associative $\Lam$-torus. Then $\aa$ admits a pre-Chevalley involution.}
\end{thm}

{As explained above to give a proof, we need to have an explicit description of structure constants of the algebra product in each case with respect to a properly chosen basis. We do it in what follows. Then Propositions \ref{pro-2}, \ref{pro2} and \ref{pro1} below all together give the proof of Theorem \ref{thm-2}.}


\subsection{Quantum tori}\label{Quantum tori}
We begin by recalling the  definition of a quantum torus from \cite[]{BGK96}.\
\begin{DEF}\label{def5}
	Let $\bq=(q_{ij})\in M_ n(\bbbk)$ be a $( n\times n)$-matrix satisfying $q_{ii}=1=q_{ij}q_{ji}$ for $1\leq i,j\leq n$. The {\it quantum torus} $\bbbk_\bq=\bbbk_\bq[x_1^{\pm1},\ldots,x_n^{\pm1}]$ based on $\bq$ is by definition the unital associative algebra with $2 n$ generators $x_1^{\pm1},\ldots,x_ n^{\pm1}$ and defining relations $x_ix_i^{-1}=1=x_i^{-1}x_i$ and $x_ix_j=q_{ij}x_jx_i$ for $1\leq i,j\leq n$. When $q_{ij}=1$ for all $i,j$, then $\bbbk_\bq=\bbbk[x_1^{\pm1},\ldots,x_n^{\pm1}]$ is just the algebra of Laurent polynomials over $\K$ in $n$-variables. Let $\{\lam_1,\ldots,\lam_n\}$ be a $\bbbz$-basis of $\Lam$. We note that $\bbbk_\bq$ is $\Lam$-graded with $\bbbk_\bq=\sum_{\a\in\Lam}\bbbk x^\a$, where for $\a=\a_1\lam_1+\cdots+\a_n\lam_n\in\Lam$,
	$x^\a:=x_1^{\a_1}\cdots x_ n^{\a_ n}$. This $\Lam$-grading makes $\K_\bq$ into an associative $\Lam$-torus. When $n=2$,
	the matrix $\bq$ is determined by a unique $q\in\bbbk^\times$ and so in this case we write $\bbbk_\bq=\bbbk_q$ with $q\in \bbbk^\times.$
\end{DEF}

\begin{pro}\label{pro-2}
	{Any associative $\Lam$-torus admits a pre-Chevalley involution.}
\end{pro}
\proof
{Any associative $\Lam$-torus is graded-isomorphic to a quantum torus $\bbbk_\bq=\sum_{\a\in\Lam}\bbbk x^\a$ for some $\bq$. Furthermore, any commutative associative $\Lam$-torus is graded-isomorphic to the algebra of Laurent polynomials over $\K$ in $n$-variables, see \cite[Lemma 1.8]{BGKN95}.
	Now, one easily sees that the assignment $x_i\mapsto x_i^{-1}$, $1\leq i\leq n$, induces a pre-Chevalley involution for $\bbbk_\bq$.}
\qed

{An {\it anti-involution} of $\aa$ is an order-2 anti-automorphism. If it also preserves the homogeneous spaces, it is called a $\Lambda$-{\it grading anti-involution}.}

\begin{rem}\label{rem1}
	Let $\bq=(q_{ij})\in M_ n(\bbbk)$ with $q_{ii}=1$ and $q_{ij}=\pm1=q_{ji}$ for $1\leq i,j\leq n$, and let $\be=(e_1,\ldots,e_n)\in\bbbk^n$ with $e_i=\pm1$ for $1\leq i\leq n$. Then there exists a unique anti-involution $\sg_\be:\bbbk_\bq\rightarrow\bbbk_\bq$ such that $\sg_\be(x_i)=e_ix_i$ for $1\leq i\leq n$, which is called the {\it anti-involution determined by} $\be$. Also $(\bbbk_\bq,\sg_\be)$ is called the {\it quantum torus with anti-involution determined by} $\be$ and $\bq$ (see \cite[p.163]{AG01}). 
\end{rem}

\subsection{Alternative tori}\label{Alternative tori}
Note that any associative $\Lam$-torus is an alternative $\Lam$-torus. But associative tori do not give all alternative tori. The octonion torus which was first found in \cite{BGKN95}, is an example of a nonassociative alternative torus. It is defined as an octonion algebra by the {\it Cayley-Dickson} process over an algebra of Laurent polynomials. Y. Yoshii described the octonion torus in a simple way via a presentation \cite{Yos08}.

\begin{exa}\label{ex1}
	Let $\cc$ be the alternative algebra over $\K$ with generators $x_i^{\pm1},1\leq i\leq
	3$, and defining relations $x_ix_i^{-1}=1=x_i^{-1}x_i$ for all $i$, $x_ix_j=-x_jx_i$ for $i\neq j$ and $(x_1x_2)x_3=-x_1(x_2x_3)$. For $n\geq3$ set $\mathbb{O}:=\cc\otimes\K[x_4^{\pm1},\ldots,x_n^{\pm1}]$ where $\K[x_4^{\pm1},\ldots,x_n^{\pm1}]$ is the algebra of Laurent polynomials over $\K$ in $n-3$ variables. Let $\epsilon_1,\ldots,\epsilon_n$ be the standard basis of $\Z^n$. Then $\mathbb{O}$ is an alternative $\Z^n$-torus, called the {\it octonion $n$-torus}, where the $\Z^n$-grading is given by $\deg(x_i)=\epsilon_i$. 
\end{exa}

\begin{rem}\label{rem3}
	The octonion $3$-torus $\cc$, also called the {\it Cayley torus}, has the following more concrete description \cite[Example 9.2]{AF11}:
	
	Let $x\in\cc$. Then one may uniquely write $x=\sum_\a k_\a x^\a$ where $x^\a=(x_1^{\a_1}x_2^{\a_2})x_3^{\a_3},\a=(\a_1,\a_2,\a_3)\in\Z^3$ and $k_\a\in\K$. Define $\epsilon:\Z^3\times\Z^3\rightarrow\Z$ by $\epsilon(\a,\b)=\a_3\b_1+\a_2\b_1+\a_3\b_2+\a_1\b_2\b_3+\a_2\b_1\b_3+\a_3\b_1\b_2$. Then the multiplication rule of $\cc$ is given by
	$$x^\a x^\b=(-1)^{\epsilon(\a,\b)}x^{\a+\b}.$$
\end{rem}

The following theorem gives the classification of alternative tori over $\K$.

\begin{thm}\cite[Theorem 1.25]{BGKN95}\cite[Corollary 5.13]{Yos02}\label{thm3}
	Any alternative $\Lam$-torus over $\K$ is graded-isomorphic to either a quantum torus or the octonion $n$-torus $\mathbb{O}$. 
\end{thm}


\begin{pro}\label{pro2}
	Let $\aa=\bigoplus_{\lam\in\Lam}\aa^\lam$ be an alternative $\Lam$-torus over $\bbbk$. Then $\aa$ admits a pre-Chevalley involution. 
\end{pro}
\proof
By Theorem \ref{thm3}, we can assume that $\aa=\bbbk_\bq$ for some $\bq$, or $\aa=\mathbb{O}$, the octonion $n$-torus. First, if $\aa=\bbbk_\bq$ then by Remark \ref{rem2} (ii) our claim is proved. Next let $\aa=\mathbb{O}=\cc\otimes\K[x_4^{\pm1},\ldots,x_n^{\pm1}]$. Then the assignment $x^\a\mapsto x^{-\a},\a\in\Z^3$ defines a pre-Chevalley  involution $\tau$ on $\cc$, since $(-1)^{\epsilon(\a,\b)}=(-1)^{\epsilon(-\a,-\b)}$, see Remark \ref{rem3}. Now we extend $\tau$ {to} $\mathbb{O}$ by $x\otimes x^\gamma\mapsto\tau(x)\otimes x^{-\gamma}$, where $x\in \cc,x^\gamma=x_4^{\gamma_4}\cdots x_n^{\gamma_n}\in\K[x_4^{\pm1},\ldots,x_n^{\pm1}]$. {It defines} a pre-Chevalley involution for $\mathbb{O}$; this completes our proof.\qed

\begin{rem}\label{rem2}
	Let $\mathbb{O}=\bigoplus_{\lam\in\bbbz^n}\mathbb{O}^\lam$ be the octonion $n$-torus. Then there exists an anti-involution $\sg$ on $\mathbb{O}$ such that $\sg(x)=\pm x$ for all $x\in\mathbb{O}^\lam$ and for all $\lam\in\bbbz^n$, see \cite[Lemma 1.20]{BGKN95}. It is called the {\it standard anti-involution} on $\mathbb{O}$.
\end{rem}

\subsection{Jordan tori}
In Example \ref{ex4} we briefly describe four families of Jordan $\Lam$-tori. In particular, we specify a basis $\{x^\a\}_{\a\in\Lam}$ and the multiplication rule in each family. For more details, see \cite{Yos00}.

\begin{DEF}\label{def3}
	Let $\aa$ be a unital commutative associative algebra, $M$ be a left $\aa$-module and $f:M\times M\rightarrow \aa$ be a symmetric $\aa$-bilinear form of $M$. Set $\mathrm{Cliff}(f):=\aa1\oplus M$. Then $\mathrm{Cliff}(f)$ together with the multiplication given by
	$$(a1,x)\cdot(b1,y)=ab1+f(x,y)1+ay+bx$$
	is a Jordan algebra. It is called the {\it Clifford Jordan algebra of} $f$ or just the {\it Jordan algebra of} $f$ (see \cite{Mc04}).
\end{DEF}

\begin{exa}\label{ex4}
	(i) First let $\bbbk_\bq=\bbbk_\bq[x_1^{\pm1},\ldots,x_n^{\pm1}]$ be the quantum torus based on a quantum matrix $\bq$ with {the} $\Lam$-grading as in Definition \ref{def5}. Now consider a new multiplication $\cdot$ on $\bbbk_\bq$ defined by $x\cdot y=\frac{1}{2}(xy+yx)$. Then $(\bbbk_\bq,\cdot)$ is a Jordan algebra denoted by $\bbbk_\bq^+$ and called the {\it plus algebra} of $\bbbk_\bq$. The $\Lam$-graded algebra $\bbbk_\bq^+=\bigoplus_{\a\in\Lam}\bbbk x^\a$ is a Jordan $\Lam$-torus and we have
	\begin{equation}\label{eq3}
		x^\a\cdot x^\b=\frac{1}{2}\prod_{i<j}q_{ij}^{\a_j\b_i}\left(1+\prod_{i,j}q_{ij}^{\a_i\b_j}\right)x^{\a+\b},
	\end{equation}
	for $\a,\b\in\Lam$ (see\cite[Example 3.2]{Yos00}).
	
	(ii) Assume that $\be=(e_{ij})$ is an elementary quantum matrix i.e., $\be$ is a quantum matrix such that $e_{ij}=1$ or $-1$ for all $i,j$. Now let $\bbbk_\be=\bbbk_\be[x_1^{\pm1},\ldots,x_n^{\pm1}]$ be the quantum torus determined by $\be$ with the unique involution $-$ such that $\overline{x_i}=x_i$ for all $i$. Then the subspace of symmetric elements $H(\bbbk_\be,-)=\{x\in\bbbk_\be\mid\overline{x}=x\}$ is a Jordan subalgebra of $\bbbk_\be^+$. In fact the $\Lam$-graded algebra $H(\bbbk_\be,-)=\bigoplus_{\a\in\Lam}\left(\bbbk x^\a\cap H(\bbbk_\be,-)\right)$ is a Jordan $\Lam$-torus (see\cite[Example 4.3 (2)]{Yos00}).
	
	(iii) Suppose that $\Lam$ is of rank $n\geq2$ and let $2\leq m\leq n$. 
	Choose free abelian subgroups $\Lambda_m$ and $\Lambda_{n-m}$ with
	$\Lambda=\Lambda_m\oplus\Lambda_{n-m}$. Let $S$ be a semilattice in $\Lam_m$. Next set $\Gamma=2\Lam_m\oplus\Lam_{n-m}$ and let $Z:=\bbbk[\Gamma]=\bigoplus_{\gamma\in\Gamma}\bbbk z^\gamma$ be the group algebra of $\Gamma$ over $\bbbk$. Now let $I$ be a set of coset representatives of $\{s+2\Lam_m\mid s\in S\}\setminus\{2\Lam_m\}$. Next let $V$ be a free $Z$-module with basis $\{t_\epsilon\}_{\epsilon\in I}$ and define a $Z$-bilinear form $f:V\times V\rightarrow Z$ by
	$$f(t_\epsilon,t_\eta)=\left\{\begin{array}{ll}
		z^{2\epsilon}&\hbox{if }\epsilon=\eta\\
		0&\hbox{otherwise}
	\end{array}\right.$$
	for all $\epsilon,\eta\in I$. Let $\jj_S:=Z\oplus V$ be the Jordan algebra over $Z$ of $f$, see Definition \ref{def3}. Note that if $\a\in S\oplus\Lam_{n-m}$ then there exist unique $\a^\prime\in\Gamma$ and $\epsilon_\a\in I\cup\{0\}$ such that $\a=\a^\prime+\epsilon_\a$. Set $t_0:=1$ and
	$$x^\a:=\left\{\begin{array}{ll}
		z^{\a^\prime}t_{\epsilon_\a}&\hbox{if }\a\in S\oplus\Lam_{n-m}\\
		0&\hbox{otherwise.}
	\end{array}\right.$$
	Then $\jj_S=\bigoplus_{\a\in\Lam}\bbbk x^\a$. In fact $\jj_S$ is a Jordan $\Lam$-torus over $\bbbk$ which  is called the {\it  standard Clifford torus determined by} $S$. For $\a,\b\in S\oplus\Lam_{n-m}$ the multiplication rule is given by (see\cite[Example 5.2)]{Yos00}:
	\begin{equation}\label{eq4}
		x^\a x^\b=\left\{\begin{array}{ll}
			x^{\a+\b}&\hbox{if }\epsilon_\a=\epsilon_b\neq0\\
			x^{\a+\b}&\hbox{if }\epsilon_\a=0\hbox{ or }\epsilon_\b=0\\
			0&\hbox{otherwise.}
		\end{array}\right.
	\end{equation}
	
	(iv) Assume that $\Lam$ is of rank $n\geq3$. Let $\omega$ be a primitive third root of unity and let $w=(w_{ij})$ be the quantum $n\times n$ matrix with $w_{12}=\omega$, $w_{21}=\omega^{-1}$ and all other entries equal to $1$. Next let $\bbbk_w=\K_w[u_1^{\pm1},\ldots,u_n^{\pm1}]$ be the quantum torus determined by $w$. Also let $Z=Z(\bbbk_w)$ be the center of $\bbbk_w$. Then $Z=\K[u_1^{\pm3},u_2^{\pm3},u_3^{\pm1},\ldots,u_n^{\pm1}]$. Now consider the first Tits construction $\mathbb{A}_t=(\bbbk_w,u_3)=\bbbk_w\oplus\bbbk_w\oplus\bbbk_w$ over $Z$, for details see \cite[Section 6]{Yos00}. Next let $\{\lam_1,\ldots,\lam_n\}$ be a basis of $\Lam$ and set $\Delta:=\bbbz\lam_1+\bbbz\lam_2+3\bbbz\lam_3+\bbbz\lam_4+\cdots+\bbbz\lam_n$. Define $\deg(u^\d=u_1^{\d_1}u_2^{\d_2}u_3^{\d_3}\cdots u_n^{\d_n})=\d$ for $\d=\d_1\lam_1+\d_2\lam_2+3\d_3\lam_3+\d_4\lam_4+\cdots+\d_n\lam_n\in\Delta$. {This equips  $\bbbk_w=\bigoplus_{\d\in\Delta}\bbbk u^\d$ with a $\Delta$-grading}. Finally for $\a=\a_1\lam_1+\cdots+\a_n\lam_n\in\Lam$, set
	$$x^\a=\left\{\begin{array}{ll}
		(u^\a,0,0)&\hbox{if }\a_3\equiv0\quad(\mod\;3)\\
		(0,u^{\a-\lam_3},0)&\hbox{if }\a_3\equiv1\quad(\mod\;3)\\
		(0,0,u^{\a+\lam_3})&\hbox{if }\a_3\equiv2\quad(\mod\;3).\\
	\end{array}\right.$$
	Then $\mathbb{A}_t=\bigoplus_{\a\in\Lam}\bbbk x^\a$ is a Jordan $\Lam$-torus over $\bbbk$ which is called {the} {\it Albert torus}. For $(x_0,x_1,x_2),(y_0,y_1,y_2)\in\mathbb{A}_t$ the multiplication is given by
	\begin{equation}\label{eq5}
		\begin{split}
			(x_0,x_1,x_2)(y_0,y_1,y_2)=&(x_0\cdot y_0+\overline{x_1y_2}+\overline{y_1x_2},\\
			&\overline{x_0}y_1+\overline{y_0}x_1+u_3^{-1}x_2\times y_2,\\
			&y_2\overline{x_0}+x_2\overline{y_0}+u_3x_1\times y_1),
		\end{split}
	\end{equation}
	where $x\cdot y=\frac{1}{2}(xy+yx),$ $x\times y=x\cdot y-\frac{1}{2}\tr(x)y-\frac{1}{2}\tr(y)x+\frac{1}{2}(\tr(x)\tr(y)-\tr(x\cdot y))1,$ $\overline{x}=x\times1$ for $x,y\in\bbbk_w$, and $\tr$ is the generic trace of the central closure $\overline{\bbbk}_w$.
\end{exa}

Y. Yoshii classified Jordan tori over any field of characteristic $\neq2$, see \cite[Theorem 2]{Yos00}. The following is the classification of Jordan tori over $\bbbk$.

\begin{thm}\cite[Corollary 7.2]{Yos00}\label{thm1}
	Any Jordan $\Lam$-torus over $\bbbk$ is graded-isomorphic to one of the four Jordan tori
	$$\bbbk_\bq^+,H(\bbbk_\be,-),\jj_S,\mathbb{A}_t.$$
\end{thm} 

\begin{pro}\label{pro1}
	Assume that $\aa=\bigoplus_{\lam\in\Lam}\aa^\lam$ is a Jordan $\Lam$-torus over $\bbbk$. Then $\aa$ admits a pre-Chevalley involution. 
\end{pro}
\proof 
By Theorem \ref{thm1}, it is enough to prove the statement for Jordan tori $\bbbk_\bq^+$, $H(\bbbk_\be,-)$, $\jj_S$, and $\mathbb{A}_t$. So let $\aa=\bigoplus_{\a\in\Lam}\bbbk x^\a$ be one of these four {families} of Jordan $\Lam$-tori and let $\{x^\a\}_{\a\in\Lam}$ be the $\bbbk$-basis of $\aa$ described in Example \ref{ex4}. We define $\tau:\aa\rightarrow\aa$ by $\bbbk$-linear extension of $\tau(x^\a)=x^{-\a}$ for all $\a\in\Lam$. Hence $\tau$ is of order $2$ and $\tau(\aa^\a)=\aa^{-\a}$ for all $\a\in\Lam$. It remains to show that $\tau$ is an automorphism of $\aa$. Let $\a,\b\in\Lam$ then $x^\a x^\b=k(\a,\b)x^{\a+\b}$ with $k(\a,\b)\in\bbbk$. For $\bbbk_\bq^+$, $H(\bbbk_\be,-)$ and $\jj_S$, it is easy to see from the multiplication rules \ref{eq3} and \ref{eq4} that $k(\a,\b)=k(-\a,-\b)$ for all $\a,\b\in\Lam$. Now let $\aa=\mathbb{A}_t$ and let $u^\a\in\K_w$. Then $u^\a=zu_1^{i_1}u_2^{i_2}$ for some $z\in Z=Z(\K_w)$. We first note that since $\tr$ is $Z$-linear then $\tr(u^\a)=z\tr(u_1^{i_1}u_2^{i_2})$. On the other hand $\tr(u_1^{i_1}u_2^{i_2})=0$ if $i_1\not\equiv0$ or $i_2\not\equiv0$ $(\mod\;3)$, see \cite[Proposition 6.7]{Yos00}. Hence
$$\tr(u^{-\a})=\left\{\begin{array}{ll}
	0&\hbox{if }i_1\not\equiv0\hbox{ or }i_2\not\equiv0\;(\mod\;3)\\
	u^{-\a}\tr(1)&\hbox{otherwise.}
\end{array}\right.$$
Considering this fact, now it is easy to check from \ref{eq5} that $k(\a,\b)=k(-\a,-\b)$ for all $\a,\b\in\Lam$. So $\tau$ is an involution of $\aa$ and this proves our claim.\qed

\section{\bf Chevalley involutions for {centerless Lie tori}}\setcounter{equation}{0}
\label{Chevalley involutions for centerless Lie tori}

{It is well known that a Lie torus is a root graded Lie algebra. Based on the structure theory of root graded Lie algebras \cite{BM92}, \cite{BZ96}, the coordinatization theorems for centerless Lie tori of reduced types have been proved in \cite{BGK96},\cite{BGKN95}, \cite{Yos00} and \cite{AG01}. Using these coordinatization theorems, we investigate the existence of Chevalley involutions for centerless Lie tori of reduced types in this section. We proceed with a {type-dependent} argument.} 

%


\begin{thm}\label{thm-3}
{Suppose $\ll$ is a centerless Lie torus of reduced type where if it is of type $A_\ell$, $\ell\geq 2$, its coordinate algebra is equipped with a $\Lambda$-grading anti-involution. Then $\ll$ admits a Chevalley involution.}
\end{thm}
The proof is provided in Propositions \ref{pron-2}, \ref{pron-6}, \ref{pro5}, \ref{pro6}, \ref{pro7}, \ref{pro8} and \ref{pro9} below.

\subsection{Simply laced centerless Lie tori of rank {$\mathbf{>2}$}}\label{simply}

Suppose in this part that $\Delta$ is the finite root system of type $A_\ell$ ($\ell\geq3$), $D_\ell$ ($\ell\geq 4$) or $E_{6,7,8}$ and let $Q=Q(\Delta)$ be the root lattice of $\Delta$. For type $A_\ell$ we let
$$\Delta=\{0\}\cup\{\epsilon_i-\epsilon_j\mid1\leq i\neq j\leq\ell+1\}$$
where $\epsilon_1,\ldots,\epsilon_{\ell+1}$ is a basis for a vector space containing $\Delta$. 

Let $\aa=\bigoplus_{\lam\in\Lam}\aa^\lam$ be an associative $\Lam$-torus and $[\aa,\aa]$ be the Lie algebra of commutators in $\aa$. Here we recall the definition of {\it special linear} Lie algebra over $\aa$.

\begin{DEF}\label{def6}	
	We denote by $M_{\ell+1}(\aa)$, the algebra of $(\ell+1)\times (\ell+1)$-matrices with entries from the associative $\Lam$-torus $\aa$, and by $\frak{gl}(\aa)$ the corresponding Lie algebra under the commutator bracket.  We set
	$$\frak{sl}_{\ell+1}(\aa):=\{A\in M_{\ell+1}(\aa)\mid\tr(A)\in [\aa,\aa]\}$$ which is a Lie subalgebra of $\frak{gl}(\aa)$ called the $(\ell+1)\times (\ell+1)$-special linear Lie algebra over $\aa$. When $\aa$ is commutative, we usually identify $\frak{sl}_{\ell+1}(\aa)$ with
	$\frak{sl}_{\ell+1}(\bbbk)\otimes_\bbbk\aa$. We note that
	$M_{\ell+1}(\aa)$, and consequently $\frak{gl}(\aa)$, inherits the $\Lam$-grading of $\aa$. In fact, we have  $\frak{gl}_{\ell+1}(\aa)=\sum_{\lam\in\Lam}\frak{gl}_{\ell+1}(\aa)^\lam$ with $\frak{gl}_{\ell+1}(\aa)^\lam=M_{\ell+1}(\aa^\lam)$. Finally, $\frak{sl}_{\ell+1}(\aa)$ is a $\Lam$-graded subalgebra of $\frak{gl}_{\ell+1}(\aa)$. 
\end{DEF}

Let $\ll=\frak{sl}_{\ell+1}(\aa)$ or $\ll=\fg\otimes\aa_{[n]}$ where $\fg$ is the finite dimensional simple Lie algebra of type $D_\ell$ ($\ell\geq 4$) or $E_{6,7,8}$. Then $\ll$ is a $Q\times\Lam$-graded Lie algebra. If $\ll=\frak{sl}_{\ell+1}(\aa)$ then
\begin{equation}\label{eq16}
	\ll_{\epsilon_i-\epsilon_j}^\lam=\aa^\lam E_{ij},\quad\ll_0^\lam=\{\Sigma_{i=1}^{\ell+1}x_iE_{ii}\mid x_i\in\aa^\lam,\Sigma_{i=1}^{\ell+1}x_i\in[\aa,\aa]\},
\end{equation}
for $i\neq j,\lam\in\Lam$, where $E_{ij}$ denotes the matrix with $1$ in {$(i,j)$}-th entry and zero elsewhere, and $\ll_\b^\lam=\{0\}$ for $\lam\in\Lam$, $\b\in Q\setminus\Delta$. Also if $\ll=\fg\otimes\aa_{[n]}$ then 
\begin{equation}
	\ll_\a^\lam=\fg_\a\otimes\aa_{[n]}^\lam,\quad\ll_\b^\lam=\{0\},
\end{equation}
for $\a\in\Delta,\b\in Q\setminus\Delta,\lam\in\Lam$.
In fact $\ll$ is a centerless $\Lam$-torus of type $\Delta$. Moreover, S. Berman, Y. Gao, and Y. Krylyuk proved the following coordinatization theorem.

\begin{thm}\cite[Theorem 2.65, Theorem 1.37]{BGK96}\label{thm0}
	(i) Any centerless Lie $\Lam$-torus of type $A_\ell$ ($\ell\geq3$) is graded-isomorphic to $\frak{sl}_{\ell+1}(\aa)$ where $\aa$ is an associative $\Lam$-torus.
	
	(ii) Any centerless Lie $\Lam$-torus of type $D_\ell$ ($\ell\geq 4$) or $E_{6,7,8}$ is graded-isomorphic to $\fg\otimes\aa_{[n]}$ where $\fg$ is the finite dimensional simple Lie algebra of type $D_\ell$ ($\ell\geq 4$) or $E_{6,7,8}$, respectively.
\end{thm}

\begin{pro}\label{pron-2}
{Let $\ll$ be a centerless Lie $\Lam$-torus of type $A_\ell$, ($\ell\geq3$),  $D_\ell$ ($\ell\geq 4$) or $E_{6,7,8}$. Further for type $A_\ell$ assume that the coordinate algebra of $\ll$ is equipped with a $\Lam$-grading anti-involution. Then $\ll$ admits a Chevalley involution.}
\end{pro}

\proof First let $\ll$ be a centerless Lie $\Lam$-torus of type  $A_\ell,$ $\ell\geq3$ {whose coordinate algebra has a $\Lam$-grading anti-involution}. From Theorem \ref{thm0}(i) we see that $\ll\cong \frak{sl}_{\ell+1}(\aa)$ where $\aa$, the coordinate algebra of $\ll$, is an associative $\Lam$-torus {with a $\Lam$-grading anti-involution $\sg$}. Thus we may assume that $\ll=\frak{sl}_{\ell+1}(\aa)$. We consider the $Q\times\Lam$-grading of $\ll=\sum_{\a\in Q,\lam\in\Lam}\ll_\a^\lam$ as in (\ref{eq16}).

Next let $\tau$ be a pre-Chevalley involution of $\aa$, see Proposition \ref{pro-2}. We define $\bar\tau:\ll\rightarrow\ll$ by $\bar\tau:=\hat\tau\circ\theta$, where $\theta$ is the involution of $\frak{sl}_{\ell+1}(\aa)$ given by
$$\theta(A)=-{\sg(A^t)},\quad A\in \frak{sl}_{\ell+1}(\aa),$$
{with $\sg(X):=(\sg(x_{ij}))$ for $X=(x_{ij})\in\mathrm{Mat}_\ell(\aa)$},
and $\hat\tau$ is the involution of $\frak{sl}_{\ell+1}(\aa)$ induced by $\tau\in\Aut(\aa)$. In fact $\hat\tau$ is given by
\begin{equation*}
	{\hat\tau(xE_{ij})=\tau(x)E_{ij},\quad x\in\aa,\;1\leq i,j\leq\ell+1.}
\end{equation*}
Then clearly $\bar\tau(\ll^\lam)=\ll^{-\lam}$ for $\lam\in\Lam$ and $\bar\tau(h)=-h$ for $h\in\ll^0_0$. So $\bar\tau$ is a Chevalley involution for $\ll={\frak{sl}_{\ell+1}(\aa)}$. 

Finally, let $\ll$ be a centerless Lie $\Lam$-torus of type $\Delta=D_\ell,\ell\geq4$ or $\Delta=E_\ell,\ell=6,7,8$. By Theorem \ref{thm0}(ii), we may assume $\ll=\fg\otimes\aa$, where $\fg$ is the finite dimensional simple Lie algebra of type $\Delta$ and $\aa=\aa_{[n]}$. Replacing $\theta$ above with $\theta_\fg\otimes1_\aa$, where $\theta_\fg$ is a Chevalley involution for $\fg$, one sees that the discussion above is valid also for this case (the involution $\hat\tau$ can be expressed as $\hat\tau(a\otimes x)=a\otimes\hat\tau(x)$ for $a\in\fg$ and $x\in \aa$).\qed

{In order to realize the compact forms of intersection matrix Lie algebras, in \cite{Gao96} the author introduces a notion of an involution for a $\Delta$-graded Lie algebra which is compatible with the $\Delta$-grading in the sense of \cite[Definition 2.4]{Gao96}. It is worth mentioning that any Chevalley involution of a Lie torus of type $\Delta$, considered as a $\Delta$-graded Lie algebra, is compatible with its $\Delta$-grading.
The involution $\theta$ is the proof of Proposition \ref{pron-2} coincides with the involution given in \cite[Example 2.5]{Gao96}, while the involution $\theta_\fg\otimes 1_\aa$ coincides with one given in \cite[Example 2.6]{Gao96}.}

\subsection{Lie tori of type $\mathbf{A_2}$}
We first recall definition of the Lie algebra $\frak{psl}_3(\aa)$ for an alternative algebra $\aa$.

Let $\aa$ be an alternative algebra and let $L_a$ and $R_a$ denote the left and the right multiplication by $a\in \aa$, respectively. For $a,b\in \aa,$ let $D_{a,b}:=[L_a,L_b]+[R_a,R_b]+[L_a,R_b]$ which is called the {\it inner derivation of $\aa$ determined by elements $a,b$}. Then, set $D_{\aa,\aa}:=\span_{\bbbk}\{D_{a,b}\mid a,b\in\aa\}$. In fact $D_{\aa,\aa}$ is an ideal of the Lie algebra $\mathrm{Der}(\aa)$ of derivations of $\aa$.

\begin{DEF}
	Assume that $\aa$ is an alternative algebra over $\K$. Set 
	$$\frak{psl}_3(\aa):=(\frak{sl}_3(\K)\otimes \aa)\oplus D_{\aa,\aa}.$$
	Then $\frak{psl}_3(\aa)$ is a Lie algebra under the bracket defined by
	\begin{equation*}
		\begin{split}
			[x\otimes a,y\otimes b]=&[x,y]\otimes\frac{ab+ba}{2}+\frac{1}{3}\tr(xy)D_{a,b}\\
			+&(xy+yx-\frac{2}{3}\tr(xy)I)\otimes\frac{ab-ba}{2},\\
			[D_{a,b},x\otimes c]=&x\otimes D_{a,b}(c)=-[x\otimes c,D_{a,b}],\\
			[D_{a,b},D_{c,d}]=&D_{D_{a,b}(c),d}+D_{c,D_{a,b}(d)},
		\end{split}
	\end{equation*}
	for $a,b,c,d\in \aa,x,y\in\frak{sl}_3(\K)$ where $I$ is the $3\times3$ identity matrix and $\tr(x)$ denotes the trace of a matrix $x$, see \cite[Section 2]{BGKN95}. 
\end{DEF}

Let $\fh$ be the Cartan subalgebra of $\frak{sl}_3(\K)$ consisting of diagonal matrices and for $i=1,2,3$, let $\epsilon_i:\fh\rightarrow\K$ be the projection onto $(i,i)$ entry. Let $\Delta:=\{0\}\cup\{\epsilon_i-\epsilon_j\mid1\leq i\neq j\leq3\}$. Then $\Delta$ can be realized as the root system of type $A_2$ with  $\{\a_1=\epsilon_1-\epsilon_2,\a_2=\epsilon_2-\epsilon_3\}$
as a basis of $\Delta$. Now let $\aa=\bigoplus_{\lam\in\Lam}\aa^\lam$ be an alternative $\Lam$-torus and set $\ll=\frak{psl}_3(\aa)$. Then $\ll$ is a $Q\times\Lam$-graded algebra with
\begin{equation}\label{eq6}
	\ll_{\epsilon_i-\epsilon_j}^\lam=E_{ij}\otimes \aa^\lam,\quad\ll_\a^\lam=\{0\},\quad\ll_0^\lam=(\fh\otimes \aa^\lam)\oplus\sum_{\mu+\nu=\lam}D_{\aa^\mu,\aa^\nu},
\end{equation}
for $i\neq j,\lam\in\Lam$ and $\a\in Q\setminus\Delta$, where $E_{ij}$ denotes the matrix with $1$ in {$(i,j)$}-th entry and zero elsewhere. In fact $\ll$ is a centerless $\Lam$-torus of type $A_2$, see \cite{BGKN95},\cite{Yos02}. The following theorem determines the coordinate algebra of EALAs of type $A_2$.

\begin{thm}\cite[Lemma 3.25]{BGKN95}\cite[Proposition 6.3]{Yos02}\label{thm4}
	Any centerless Lie $\Lam$-torus of type $A_2$ is graded-isomorphic to $\frak{psl}_3(\aa)$ where $\aa$ is an alternative $\Lam$-torus.
\end{thm}

\begin{pro}\label{pron-6}
{Suppose $\ll$ is a centerless $\Lambda$-torus of type $A_2$ and that its coordinate algebra is equipped with a $\Lambda$-grading anti-involution. Then $\ll$ admits a Chevalley involution.}\end{pro}
\proof
Using Theorem \ref{thm4}, it is enough to show the existence of a Chevalley involution for $\ll=\frak{psl}_3(\aa)$ where $\aa=\bigoplus_{\lam\in\Lam}\aa^{\lam}$ is an alternative $\Lam$-torus {with a $\Lam$-grading anti-involution $\bar{}$}. 
{First recall from Proposition \ref{pro2} that} $\aa$ admits a pre-Chevalley involution $\tau$. Set $\bar\tau:=\tau\circ\bar{}$. Then $\bar\tau$ is an anti-automorphism of $\aa$. So $\bar\tau$ induces a Lie algebra anti-automorphism $\hat{\tau}:\mathrm{Der}(\aa)\rightarrow\mathrm{Der}(\aa)$ by $\hat{\tau}(D)=\bar\tau D\bar\tau^{-1}$, for $D\in\mathrm{Der}(\aa)$. Since $\hat{\tau}(D_{a,b})=D_{\tau(\overline{a}),\tau(\overline{b})}$ for $a,b\in \aa$, we see that $\hat{\tau}$ restricts to a Lie algebra anti-automorphism of $D_{\aa,\aa}$. Next, we define $\tilde\tau$ on $\frak{psl}_3(\aa)=(\frak{sl}_3(\K)\otimes \aa)\oplus D_{\aa,\aa}$ by
$$\tilde\tau(x\otimes a+D_{b,c})=-x^t\otimes\tau(\overline{a})+\hat{\tau}(D_{b,c}),$$
for $x\in\frak{sl}_3(\K),a,b,c\in \aa$, where $x^t$ denotes the transpose of $x$. Our claim is that $\tilde\tau$ is a Chevalley involution for $\frak{psl}_3(\aa)$. To see this, let $x,y\in\frak{sl}_3(\K),$ and $a,b\in \aa$. Then
\begin{equation*}
	\begin{split}
		\tilde\tau[x\otimes a,y\otimes b]=&-[x,y]^t\otimes\frac{\tau(\overline{ab+ba})}{2}+\frac{1}{3}\tr(xy)\hat{\tau}(D_{a,b})\\
		-&(xy+yx-\frac{2}{3}\tr(xy)I)^t\otimes\frac{\tau(\overline{ab-ba})}{2}.
	\end{split}
\end{equation*}
Now we have 
\begin{equation*}
	\begin{split}
		-[x,y]^t\otimes\frac{\tau(\overline{ab+ba})}{2}=&[x^t,y^t]\otimes\frac{\tau(\overline{a})\tau(\overline{b})+\tau(\overline{b})\tau(\overline{a})}{2},\\
		\tr(xy)\hat{\tau}(D_{a,b})=&\tr(x^ty^t)D_{\tau(\overline{a}),\tau(\overline{b})},
	\end{split}
\end{equation*}
and
\begin{equation*}
	\begin{split}
		&-(xy+yx-\frac{2}{3}\tr(xy)I)^t\otimes\frac{\tau(\overline{ab-ba})}{2}=\\&(x^ty^t+y^tx^t-\frac{2}{3}\tr(x^ty^t)I)\otimes\frac{\tau(\overline{a})\tau(\overline{b})-\tau(\overline{b})\tau(\overline{a})}{2}.
	\end{split}
\end{equation*}
So $\tilde\tau[x\otimes a,y\otimes b]=[\tilde\tau(x\otimes a),\tilde\tau(y\otimes b)]$. Since $$(\tau(\overline{b}),\tau(\overline{c}),\tau(\overline{a}))=-(\tau(\overline{a}),\tau(\overline{c}),\tau(\overline{b})),\;a,b,c\in \aa$$ 
where $(\tau(\overline{b}),\tau(\overline{c}),\tau(\overline{a}))$ is the associator of $\tau(\overline{b}),\tau(\overline{c}),\tau(\overline{a})$, we get $\tau(\overline{D_{a,b}(c)})=D_{\tau(\overline{a}),\tau(\overline{b})}(\tau(\overline{c}))$. This gives $\tilde\tau[D_{a,b},x\otimes c]=[\tilde\tau(D_{a,b}),\tilde\tau(x\otimes c)]$ and $\tilde\tau[D_{a,b},D_{c,d}]=[\tilde\tau(D_{a,b}),\tilde\tau(D_{c,d})]$. So $\tilde\tau$ is a Lie algebra involution of $\ll=\frak{psl}_3(\aa)$.

To see that $\tilde\tau$ is a pre-Chevalley involution, {fix $\lam\in\Lam$,} then by (\ref{eq6}),
$$\ll^\lam=(\fh\otimes \aa^\lam)\oplus\sum_{\mu+\nu=\lam}D_{\aa^\mu,\aa^\nu}\oplus\left(\bigoplus_{\epsilon_i-\epsilon_j\in\Delta} E_{ij}\otimes \aa^\lam\right).$$
So $\tilde\tau(\ll^\lam)=\ll^{-\lam}$. Also we note that $\ll_0^0=(\fh\otimes \K1)\oplus\sum_{\mu\in\Lam}\K D_{x^\mu,\tau(x^\mu)}$, where $0\neq x^\mu\in \aa^\mu$. Since $D_{a,b}+D_{b,a}=0$ (see \cite[Section 1]{BGKN95}) and $D_{a,b}=D_{-a,-b}$ for all $a,b\in \aa$, we have $\tilde\tau(D_{x^\mu,\tau(x^\mu)})=D_{\tau(x^\mu),x^\mu}=-D_{x^\mu,\tau(x^\mu)}$. This shows that $\tilde\tau_{|_{\ll_0^0}}=-\id$. So $\tilde\tau$ is a Chevalley involution for $\frak{psl}_3(\aa)$ and this completes the proof.\qed

\subsection{Centerless Lie tori of type $\mathbf{A_1}$ and $\mathbf{C_\ell (\ell\geq2)}$}
We begin by recalling from {\cite[Chapter VIII, Section 5]{Jac68}} the Tits-Kantor-Koecher ($\mathrm{TKK}$) construction of a Lie algebra from a Jordan algebra. We apply the $\mathrm{TKK}$ construction to a $\Lam$-graded Jordan algebra, to construct centreless $\Lam$-tori of type $A_1$ and $C_l$. 

\parag
($\mathrm{TKK}$ Construction)
Let $\jj$ be a Jordan algebra and let $L_x$ be the operator on $\jj$ defined by $L_xy=xy$, for $x,y\in\jj$. Set 
$$\mathrm{Instrl}(\jj):=L_\jj\oplus\mathrm{Inder}(\jj)$$ 
where $\mathrm{Inder}(\jj)$ is the Lie algebra of inner derivations $\{\sum_i[L_{x_i},L_{y_i}]\mid x_i,y_i\in\jj\}$. Together with the bracket defined by
$$[L_x+C,L_y+D]=[L_x,L_y]+L_{Cy}-L_{Dx}+[C,D],$$
for $x,y\in\jj,C,D\in\mathrm{Inder}(\jj)$, $\mathrm{Instrl}(\jj)$ is a subalgebra of the Lie algebra $\gl(\jj)$ which is called the {\it inner structure Lie algebra} of $\jj$. Define $\bar{}:\mathrm{Instrl}(\jj)\rightarrow\mathrm{Instrl}(\jj)$ by $\overline{L_a+D}=-L_a+D$. It is an involution of $\mathrm{Instrl}(\jj)$. Now put $\mathrm{TKK}(\jj):=\jj\oplus\mathrm{Instrl}(\jj)\oplus\bar{\jj}$ where $\bar{\jj}$ is isomorphic to $\jj$ under a linear map $x\mapsto\bar{x}$. Then $\mathrm{TKK}(\jj)$ is a Lie algebra, called {\it TKK Lie algebra} of $\jj$, under the bracket defined by
\begin{eqnarray*}
	[x_1+\bar{y}_1+E_1,x_2+\bar{y}_2+E_2]&=&-E_2x_1+E_1x_2-\overline{\bar E_2b_1}+\overline{\bar E_1b_2}\\
	&+&x_1\triangle y_2-x_2\triangle y_1+[E_1,E_2],
\end{eqnarray*}
for $x_i\in\jj,\bar y_i\in\bar{\jj}$, and $E_i\in\mathrm{Instrl}(\jj)$, where $x\triangle y=L_{xy}+[L_x,L_y]$.

Let $\jj=\bigoplus_{\lam\in\Lam}\jj^\lam$ be a $\Lam$-graded Jordan algebra. Then $\mathrm{TKK}(\jj)$ is $\Lam$-graded with
\begin{equation}\label{eq9}
	\mathrm{TKK}(\jj)^\lam=\jj^\lam\oplus\mathrm{Instrl}(\jj)^\lam\oplus\overline{\jj^\lam},
\end{equation}
for $\lam\in\Lam$, where $\mathrm{Instrl}(\jj)^\lam=L_{\jj^\lam}\oplus\sum_{\mu+\nu=\lam}[L_{\jj^\mu},L_{\jj^\nu}]$ (for more details see \cite[Chapter III, $\S 2$]{AABGP97}).

\begin{pro}\label{pron-5}
	Assume that $\jj=\bigoplus_{\lam\in\Lam}\jj^\lam$ is a $\Lam$-graded Jordan algebra with a pre-Chevalley involution. Then $\mathrm{TKK}(\jj)$ admits a pre-Chevalley involution.	
\end{pro}
\proof 
Let $\tau:\jj\rightarrow\jj$ be a pre-Chevalley involution on $\jj$. We first note that $\tau$ induces a Lie algebra involution $\hat{\tau}$ of $\mathrm{Der}(\jj)$ by $\hat{\tau}(D)=\tau D\tau^{-1}$, for $D\in\mathrm{Der}(\jj)$. Now we have
$$\hat{\tau}(\sum_i[L_{x_i},L_{y_i}])=\sum_i[L_{\tau(x_i)},L_{\tau(y_i)}],$$
where $\;x_i,y_i\in\jj$. So $\hat{\tau}$ restricts to a Lie algebra involution of $\mathrm{Inder}(\jj)$. Next, we extend $\hat{\tau}$ to $\mathrm{Instrl}(\jj)$ by $\hat{\tau}(L_x+D):=L_{\tau(x)}+\hat{\tau}(D)$, for $x\in\jj$ and $D\in\mathrm{Inder}(\jj)$. Then $\hat{\tau}$ is a Lie algebra involution of $\mathrm{Instrl}(\jj)$, since $\hat{\tau}(L_{Ex})=L_{\hat{\tau}(E)(\tau(x))},x\in\jj,E\in\mathrm{Inder}(\jj)$. Now we define $\bar{\tau}$ on $\mathrm{TKK}(\jj)= J\oplus\mathrm{Instrl}(\jj)\oplus\bar{\jj}$ by
\begin{equation}\label{eq11}
	\bar{\tau}(x+\bar{y}+E)=\overline{\tau(x)}+\tau(y)+\overline{\hat\tau(E)},\;x\in\jj,\bar{y}\in\bar{\jj},E\in\mathrm{Instrl}(\jj).
\end{equation}
To show that $\bar\tau$ is a Lie algebra homomorphism of $\mathrm{TKK}(\jj)$, let $x_1,x_2\in\jj,\bar{y}_1,\bar{y}_2\in\bar{\jj},E_1,E_2\in\mathrm{Instrl}(\jj)$. Then
\begin{eqnarray*}
	\bar\tau[x_1+\bar{y}_1+E_1,x_2+\bar{y}_2+E_2]&=&-\overline{\tau(E_2x_1)}+\overline{\tau(E_1x_2)}-\tau(\bar E_2y_1)+\tau(\bar E_1y_2)\\
	&+&\overline{\hat\tau(x_1\triangle y_2)}-\overline{\hat\tau(x_2\triangle y_1)}+[\overline{\hat\tau(E_1)},\overline{\hat\tau(E_2)}].
\end{eqnarray*}
We have $\overline{\tau(E_2x_1)}=\overline{\hat\tau(E_2)\tau(x_1)}$. Similarly $\overline{\tau(E_1x_2)}=\overline{\hat\tau(E_1)\tau(x_2)},\tau(\bar E_2y_1)=\overline{\hat\tau(E_2)}\tau(y_1)$ and $\tau(\bar E_1y_2)=\overline{\hat\tau(E_1)}\tau(y_2)$. On the other hand
\begin{eqnarray*}
	\overline{\hat\tau(x_1\triangle y_2)}&=&-L_{\tau(x_1y_2)}+[L_{\tau(x_1)},L_{\tau(y_2)}]\\
	&=&-(L_{\tau(y_2x_1)}+[L_{\tau(y_2)},L_{\tau(x_1)}])=-\tau(y_2)\triangle\tau(x_1).
\end{eqnarray*}
Similarly $\overline{\hat\tau(x_2\triangle y_1)}=-\tau(y_1)\triangle\tau(x_2)$. So $\bar\tau$ is an involution of $\mathrm{TKK}(\jj)$. Finally, note that $\bar\tau(\mathrm{TKK}(\jj)^\lam)=\mathrm{TKK}(\jj)^{-\lam}$ for all $\lam\in\Lam$, see (\ref{eq9}). Thus $\bar\tau$ is a Chevalley involution for $\ll=\mathrm{TKK}(\jj)$ and this proves our claim.
\qed

\subsubsection{\bf Type $\mathbf{A_1}$}

Suppose that $\Delta$ is the finite root system of type $A_1$ i.e. $\Delta=\{0,\pm\a\}$ and let $Q=\Z\a$.

Next assume that $\jj=\bigoplus_{\lam\in\Lam}\jj^\lam$ is a Jordan $\Lam$-torus, see Definition \ref{def7}. Now set $\ll=\mathrm{TKK}(\jj)$ and let
\begin{equation}\label{eq2}
	\ll_\a=\jj,\quad\ll_{-\a}=\bar{\jj},\quad\ll_\b=\{0\},\quad\ll_0=\mathrm{Instrl}(\jj),
\end{equation}
for $\b\in Q\setminus\Delta$. Considering the $\Lam$-grading (\ref{eq9}), then $\ll$ is $Q\times\Lam$-graded. Y. Yoshii showed that $\ll=\mathrm{TKK}(\jj)$ is a centreless Lie $\Lam$-torus of type $A_1$, see \cite{Yos00}. Moreover he proved the following result in which he determined the coordinate algebra of EALAs of type $A_1$.

\begin{thm}\cite[Theorem 1]{Yos00}\label{thm2}
	Any centreless Lie $\Lam$-torus of type $A_1$ is graded-isomorphic to $\mathrm{TKK}(\jj)$ where $\jj$ is a Jordan $\Lam$-torus.
\end{thm}

\begin{pro}\label{pro5}
	Let $\ll$ be a centerless $\Lam$-torus of type $A_1$. Then $\ll$ admits a Chevalley involution.	
\end{pro}
\proof
By Theorem \ref{thm2}, we can assume that $\ll=\mathrm{TKK}(\jj)$ where $\jj$ is a Jordan $\Lam$-torus. We now show that $\mathrm{TKK}(\jj)$ admits a Chevalley involution. First, we mention that there exists a pre-Chevalley involution $\tau:\jj\rightarrow \jj$, see Proposition \ref{pro1}. Then, by Proposition \ref{pron-5}, $\mathrm{TKK}(\jj)$ admits a pre-Chevalley involution $\bar{\tau}$ defined by (\ref{eq11}). Also note that $\ll_0^0=\K L_1\oplus\sum_{\mu\in\Lam}\K[L_{x^\mu},L_{\tau(x^\mu)}]\subseteq\mathrm{Instrl}(\jj)$ where $0\neq x^\mu\in\jj^\mu$, see (\ref{eq2}), (\ref{eq9}). Since $\tau$ is of order $2$, we get $\bar\tau([L_{x^\mu},L_{\tau(x^\mu)}])=-[L_{x^\mu},L_{\tau(x^\mu)}]$. On the other hand $\bar\tau(L_1)=\overline{\hat{\tau}(L_1)}=-L_1$. Hence we have $\bar\tau_{|_{\ll_0^0}}=-\id$. So $\bar{\tau}$ is a Chevalley involution for $\ll=\mathrm{TKK}(\jj)$ and this completes our proof.
\qed

\subsubsection{\bf{Type} $\mathbf{C_\ell (\ell\geq2)}$ }

Assume that 
$$\Delta=\{0\}\cup\{\pm2\ep_i\mid1\leq i\leq l\}\cup\{\pm(\ep_i\pm\ep_j)\mid1\leq i< j\leq \ell\}$$
is the finite root system of type $C_\ell$ and let $Q=Q(\Delta)$ be the root lattice of $\Delta$.

Let $\aa$ be an algebra over $\bbbk$ and let $\sg$ be an anti-involution of $\aa$. For $\ell\geq1$, set
$$H_\ell(\aa,\sg):=\{X\in\mathrm{Mat}_\ell(\aa)\mid\sg(X^t)=X\}$$
where  $\sg(X):=(\sg(x_{ij}))$ for $X=(x_{ij})\in\mathrm{Mat}_\ell(\aa)$. Then $H_l(\aa,\sg)$ is a subalgebra of $\mathrm{Mat}_\ell(\aa)^+$, the plus algebra of $\mathrm{Mat}_\ell(\aa)$. So $H_l(\aa,\sg)$ is an algebra with multiplication $X\cdot Y=\frac{1}{2}(XY+YX)$. Note that if $\aa$ is $\Lam$-graded, then the $\Lam$-grading of $\aa$ induces a $\Lam$-grading on $H_\ell(\aa,\sg)$ such that
$$H_\ell(\aa,\sg)^\lam=\{X=(x_{ij})\in\mathrm{Mat}_\ell(\aa)\mid x_{ij}\in\aa^\lam,\text{for all }1\leq i,j\leq \ell\},$$
for $\lam\in\Lam$.

\begin{lem}\label{lem2}
	(i) Let $\aa$ be an algebra with an anti-involution $\sg$ and let $\tau$ be an automorphism of $\aa$. If $\sg\tau=\tau\sg$, then $\tau$ induces an automorphism $\tau:H_\ell(\aa,\sg)\rightarrow H_\ell(\aa,\sg)$ by $\tau(X):=(\tau(x_{ij}))$ for $X=(x_{ij})\in H_\ell(\aa,\sg)$. 
	
	(ii) Consider the quantum torus $(\aa,\sg)=(\bbbk_\bq,\sg_\be)$ with anti-involution determined by some $\be$ and $\bq$ satisfying the conditions of Remark \ref{rem1}, or consider the octonion torus $(\aa,\sg)=(\mathbb{O},\sg)$ with standard anti-involution. Then $H_\ell(\aa,\sg)$ admits a pre-Chevalley involution.
\end{lem}
\proof
(i) Clearly {$\tau:\mathrm{Mat}_\ell(\aa)\rightarrow\mathrm{Mat}_\ell(\aa)$ with $\tau(X):=(\tau(x_{ij}))$ is an automorphism of $(\mathrm{Mat}_\ell(\aa),\cdot)$. Since $\sg\tau=\tau\sg$, we have $\sg(\tau(X)^t)=\sg(\tau(X^t))=\tau(\sg(X^t))=\tau(X)$. So $\tau:\mathrm{Mat}_\ell(\aa)\rightarrow\mathrm{Mat}_\ell(\aa)$ restricts to an automorphism of $H_\ell(\aa,\sg)$.}

(ii) By Theorem \ref{thm-2}, $\aa$ admits a pre-Chevalley involution $\tau$. Considering Remark \ref{rem1} and Remark \ref{rem2}, one can easily check that $\sg\tau=\tau\sg$. Then by part (i), $\tau:H_\ell(\aa,\sg)\rightarrow H_\ell(\aa,\sg)$ is an automorphism. Also we have $\tau(H_l(\aa,\sg)^\lam)=H_\ell(\aa,\sg)^{-\lam}$, for all $\lam\in\Lam$. So $\tau:H_\ell(\aa,\sg)\rightarrow H_\ell(\aa,\sg)$ is a pre-Chevalley involution. 
\qed

Let $\aa$ be a commutative associative algebra, $M$ be a left $\aa$-module and $f:M\times M\rightarrow \aa$ be a symmetric $\aa$-bilinear form of $M$. Set $\mathrm{RedCliff}(f):=\aa\oplus \aa\oplus M$. Then $\mathrm{RedCliff}(f)$ together with the multiplication given by
{\small
\begin{equation*}
	\begin{split}	
(a_1,a_2,x)\cdot(b_1,b_2,y)=\left(a_1b_1+f(x,y),a_2b_2+f(x,y),\frac{1}{2}(a_1y+a_2y+b_1x+b_2x)\right)
\end{split}
\end{equation*}
}
is a Jordan algebra. It is called the {\it reduced Clifford Jordan algebra of} $f$ (see \cite{Mc04}).

{Now assume that $\aa_{[n]}$ is the ring of Laurent polynomials in $n$ variables and $\aa_{[n]}^m$ is the free $\aa_{[n]}$-module with basis $v_1,\ldots,v_m$.} Fix $\tau_1,\ldots,\tau_m\in\bbbz^n$, $m\geq2$, such that $\tau_1=0$ and $\tau_i\not\equiv\tau_j (\mod\;2\bbbz^n)$ for $i\neq j$. Define $h:\aa_{[n]}^m\times \aa_{[n]}^m\rightarrow \aa_{[n]}$ by
\begin{equation}\label{eq8}
	h\left(\sum_{i=1}^{m}a_iv_i,\sum_{i=1}^{m}b_iv_i\right)=\sum_{i=1}^{m}a_ib_ix^{\tau_i},
\end{equation}
Then $\mathrm{RedCliff}(h)$ has the following $\Lam$-grading:
\begin{equation}
	\deg(x^\lam,0,0)=\deg(0,x^\lam,0)=2\lam
	\quadd\text{and}\quadd
	\deg(0,0,x^\lam v_i)=2\lam+\tau_i,
\end{equation}
where $\lam\in\bbbz^n$ and $1\leq i\leq m$ (see \cite[]{AG01}).

\begin{lem}\label{lem3}
	Let $h:\aa_{[n]}^m\times \aa_{[n]}^m\rightarrow \aa_{[n]}$ be as in (\ref{eq8}). Then $\mathrm{RedCliff}(h)$ admits a pre-Chevalley involution. 
\end{lem} 
\proof
Let $\tau:\mathrm{RedCliff}(h)\rightarrow\mathrm{RedCliff}(h)$ be the linear map defined by
$$(x^\lam,0,0)\mapsto(x^{-\lam},0,0),\;\;(0,x^\lam,0)\mapsto(0,x^{-\lam},0),\;\;(0,0,x^\lam v_i)\mapsto(0,0,x^{-\lam-\tau_i} v_i),$$
for $\lam\in\bbbz^n$ and $1\leq i\leq m$. Clearly $\tau(\mathrm{RedCliff}(h)^\lam)=\mathrm{RedCliff}(h)^{-\lam}$, for all $\lam\in\bbbz^n$. One checks that $\tau$ is an automorphism of $\mathrm{RedCliff}(h)$. Hence $\tau$ is a pre-Chevalley involution of $\mathrm{RedCliff}(h)$.
\qed

Let $\jj=H_\ell(\aa,\sg)$, where $(\aa,\sg)$ is as in Lemma \ref{lem2}(ii), or let $\jj=\mathrm{RedCliff}(h)$, where $h:\aa_{[n]}^m\times \aa_{[n]}^m\rightarrow \aa_{[n]}$ is given by (\ref{eq8}). In each case, $\jj$ has a set of orthogonal idempotents $\mathcal{I}=\{e_i\mid 1\leq i\leq \ell\}$ which sum to $1$. The elements of $\mathcal{I}$ for $H_\ell(\aa,\sg)$ are diagonal matrix units, and for  $\mathrm{RedCliff}(h)$ are $(1,0,0),(0,1,0)$. Thus, $\mathcal{I}$ determines a Peirce decomposition $\jj=\bigoplus_{1\leq i\leq j\leq \ell}\jj_{ij}$, where
$$\jj_{ii}=\{x\in\jj\mid e_i\cdot x=x\},\quad 1\leq i\leq \ell,$$
and
$$\jj_{ij}=\{x\in\jj\mid e_i\cdot x=e_j\cdot x=\frac{1}{2}x\},\quad 1\leq i\neq j\leq \ell.$$
This decomposition leads to a $Q$-grading of $\ll=\mathrm{TKK}(\jj)$ such that
\begin{equation}\label{eq10}
	\ll_{\ep_i+\ep_j}=\jj_{ij},\quad\ll_{-\ep_i-\ep_j}=\overline{\jj_{ij}},\quad 1\leq i\leq j\leq \ell.
\end{equation}
Considering the $\Lam$-grading (\ref{eq9}), then $\ll$ is $Q\times\Lam$-graded (\cite[Remark 4.88]{AG01}). In fact, $\ll=\mathrm{TKK}(\jj)$ is a centerless Lie $\Lam$-torus of type $C_\ell (\ell\geq2)$, see \cite[Remark 4.89]{AG01}. Moreover B. Allison and Y. Gao proved the following coordinatization theorem for type $C_\ell (\ell\geq2)$.

\begin{thm}\cite[Theorem 4.87]{AG01}\label{thm5}
	Any centerless Lie $\Lam$-torus of type $C_\ell, \ell\geq2$ is graded-isomorphic to $\mathrm{TKK}(\jj)$ where one of the following hold:
	
	(i) $\jj=H_\ell(\bbbk_\bq,\sg_\be)$, where $(\bbbk_\bq,\sg_\be)$ is the quantum torus with anti-involution determined by some $\be$ and $\bq$ satisfying the conditions of Remark \ref{rem1}.
	
	(ii) $\ell=3, n\geq3$ and $\jj=H_\ell(\mathbb{O},\sg)$, where $(\mathbb{O},\sg)$ is the octonion torus with standard anti-involution.
	
	(iii) $\ell=2$ and $\jj=\mathrm{RedCliff}(h)$, where $h:\aa_{[n]}^m\times \aa_{[n]}^m\rightarrow \aa_{[n]}$ is defined by \ref{eq8}.
\end{thm}

\begin{pro}\label{pro6}
	Assume that $\ll$ is a centerless $\Lam$-torus of type $C_\ell, \ell\geq2$. Then $\ll$ admits a Chevalley involution.	
\end{pro}
\proof
We know that $\ll$ is graded-isomorphic to $\mathrm{TKK}(\jj)$ with $\jj$ as in part (i), (ii) or (iii) of Theorem \ref{thm5}. In each case, $\jj$ has a pre-Chevalley involution, see Lemmas \ref{lem2}(ii) and \ref{lem3}. Then, by Proposition \ref{pron-5}, $\mathrm{TKK}(\jj)$ admits a pre-Chevalley involution $\bar{\tau}$ defined by (\ref{eq11}). So to justify our claim, it is enough to show that $\bar{\tau}$ acts on $\mathrm{TKK}(\jj)_0^0=\sum_{\a\in\Delta^\times,\mu\in\Lam} [\mathrm{TKK}(\jj)_\a^\mu,\mathrm{TKK}(\jj)_{-\a}^{-\mu}]\subset\mathrm{Instrl}(\jj)$ as $-\id$. Note that $\bar\tau(\mathrm{TKK}(\jj)_\a)=\mathrm{TKK}(\jj)_{-\a}$, for $\a\in\Delta^\times$, since $\bar\tau(\jj_{ij})=\overline{\jj_{ij}}$, for $1\leq i\leq j\leq \ell$, see (\ref{eq10}). Thus we have $\mathrm{TKK}(\jj)_0^0=\sum_{\a\in\Delta^\times,\mu\in\Lam}\K[x^\mu_\a,\bar\tau(x^\mu_\a)]$, where $0\neq x^\mu_\a\in\mathrm{TKK}(\jj)^\mu_\a$. Since $\bar\tau$ is of order $2$, we get $\bar\tau([x^\mu_\a,\bar\tau(x^\mu_\a)])=-[x^\mu_\a,\bar\tau(x^\mu_\a)]$. So $\bar\tau_{|_{\mathrm{TKK}(\jj)_0^0}}=-\id$ and this completes our proof.
\qed

\subsection{Centerless Lie tori of type $\mathbf{B_\ell(\ell\geq3)}$, $\mathbf{F_4}$ and $\mathbf{G_2}$}

The generalized Tits construction can be used to construct centerless Lie tori of type $B_\ell(\ell\geq3)$, $F_4$ and $G_2$. Here, we begin by recalling the generalized Tits construction (see \cite[pp.120–127]{Sch66} and \cite[$\S10$]{Jac71} for some references of this construction). Then we separate into three parts, one for each of these types.


Assume that $\aa$ is a unital commutative associative algebra over $\K$. Let $\mathcal{X}$ be a unital algebra over $\aa$. An $\aa$-linear map $T:\mathcal{X}\rightarrow \aa$ is called a {\it normalized trace} on $\mathcal{X}$ if $T(1)=1$ and
$$T(xx')=T(x'x),\quadd T((xx')x'')=T(x(x'x''))$$
for $x,x',x''\in\mathcal{X}$. We
use the same symbol $T$ to denote the normalized trace for different algebras. If $T$ is a normalized trace on $\mathcal{X}$, then one can write
$$\mathcal{X}=\aa1\oplus\mathcal{X}_0$$
where $\mathcal{X}_0=\{x\in\mathcal{X}\mid T(x)=0\}$. Also for $x,x'\in\mathcal{X}$, we have 
$$xx'=T(xx')1+x*x'$$
where $x*x'$ denotes the projection of $xx'$ onto $\mathcal{X}_0$. In the following we write $\mathrm{Der}_A^0(\mathcal{X})$ for the Lie subalgebra of the $\aa$-derivations of $\mathcal{X}$ which send $\mathcal{X}_0$ to $\mathcal{X}_0$.

\parag
(Generalized Tits Construction)
Let $\aa,\bb$ be unital commutative associative algebras over $\K$ and let $\mathcal{X},\mathcal{Y}$ be unital algebras with normalized traces over $\aa$ and $\bb$, respectively. Also let $\mathcal{D}(\mathcal{X}),\mathcal{D}(\mathcal{Y})$ be subalgebras of $\mathrm{Der}_A^0(\mathcal{X})$ and $\mathrm{Der}_B^0(\mathcal{Y})$, respectively and assume there are $\aa$ {and $\bb$}-bilinear skew-symmetric maps $(x,x^\prime)\rightarrow D_{x,x^\prime}$ and $(y,y^\prime)\rightarrow D_{y,y^\prime}$ of $\mathcal{X}_0\times\mathcal{X}_0$ into $\mathcal{D}(\mathcal{X})$ and $\mathcal{Y}_0\times\mathcal{Y}_0$ into $\mathcal{D}(\mathcal{Y})$, respectively, satisfying
$$[D,D_{x,x^\prime}]=D_{Dx,x^\prime}+D_{x,Dx^\prime},\quadd[E,D_{y,y^\prime}]=D_{Ey,y^\prime}+D_{y,Ey^\prime}$$
for $D\in\mathcal{D}(\mathcal{X}),x,x^\prime\in\mathcal{X}_0,E\in\mathcal{D}(\mathcal{Y}),y,y^\prime\in\mathcal{Y}_0$. Set
$$\mathcal{T}(\mathcal{X}/\aa,\mathcal{Y}/\bb):=(\mathcal{D}(\mathcal{X})\otimes \bb)\oplus(\mathcal{X}_0\otimes\mathcal{Y}_0)\oplus(\aa\otimes\mathcal{D}(\mathcal{Y}))$$
as a $\K$-vector space. Then $\mathcal{T}(\mathcal{X}/A,\mathcal{Y}/\bb)$ is an algebra over $\K$ with the anticommutative multiplication given by
\begin{eqnarray*}
	&[D\otimes b,D^\prime\otimes b^\prime]=[D,D^\prime]\otimes bb^\prime,\\
	&[a\otimes E,a^\prime\otimes E^\prime]=aa^\prime\otimes[E,E^\prime],\\
	&[D\otimes b,a\otimes E]=0,\\
	&[D\otimes b,x\otimes y]=Dx\otimes by=-[x\otimes y,D\otimes b],\\
	&[a\otimes E,x\otimes y]=ax\otimes Ey=-[x\otimes y,a\otimes E],\\
	&[x\otimes y,x^\prime\otimes y\prime]=D_{x,x^\prime}\otimes T(yy^\prime)+(x*x^\prime)\otimes(y*y^\prime)+T(xx^\prime)\otimes D_{y,y^\prime}
\end{eqnarray*}
for $D,D^\prime\in\mathcal{D}(\mathcal{X}),b,b^\prime\in \bb,a,a^\prime\in \aa,E,E^\prime\in\mathcal{D}(\mathcal{Y}),x,x^\prime\in\mathcal{X}_0,y,y^\prime\in\mathcal{Y}_0$. If $\mathcal{X},\mathcal{Y}$ are chosen appropriately, then $\mathcal{T}(\mathcal{X}/\aa,\mathcal{Y}/\bb)$ will be a Lie algebra, see \cite[Proposition 3.9]{BZ96}.

\subsubsection{\bf{Type} $\mathbf{B_\ell (\ell\geq3)}$ }

Assume that 
$$\Delta=\{0\}\cup\{\pm\ep_i\mid1\leq i\leq \ell\}\cup\{\pm(\ep_i\pm\ep_j)\mid1\leq i< j\leq \ell\}$$
is the finite root system of type $B_\ell$ with $\Delta_{sh}=\{\pm\ep_i\mid1\leq i\leq \ell\}$, and $\Delta_{lg}=\{\pm(\ep_i\pm\ep_j)\mid1\leq i< j\leq \ell\}$ as short and long roots, respectively and let $Q=Q(\Delta)$ be the root lattice of $\Delta$.

Let $\mathbb{V}$ be a $(2\ell+1)$-dimensional vector space over $\K$ with basis $v_1,\ldots,v_{2\ell+1}$ and let $f$ be the symmetric bilinear form on $\mathbb{V}$ defined by
\begin{equation*}
	f(v_i,v_j)=\left\{\begin{array}{ll}
		1&\hbox{if }1\leq i\leq\ell,j=\ell+i\\
		1&\hbox{if }i=j=2\ell+1\\
		0&\hbox{otherwise.}
	\end{array}\right.
\end{equation*}
Now let $\fg$ be the Lie algebra of endomorphisms of $\mathbb{V}$ which are skew relative to $f$. Then $\fg$ is the finite dimensional simple Lie algebra of type $B_\ell$ with the Cartan subalgebra $\fh$ consisting of elements of $\fg$ which are diagonal relative to the basis $v_1,\ldots,v_{2\ell+1}$. The set of weights of $\fh$-module $\mathbb{V}$ is 
$$\{0\}\cup\Delta_{sh}=\{0\}\cup\{\pm\ep_i\mid1\leq i\leq \ell\}$$ 
and its weight spaces are
\begin{equation*}
	\mathbb{V}_{\a}=\left\{\begin{array}{ll}
		\K v_i&\hbox{if }\a=\ep_i\\
		\K v_{\ell+i}&\hbox{if }\a=-\ep_i\\
		\K v_{2\ell+1}&\hbox{if }\a=0.
	\end{array}\right.
\end{equation*}

Next let $\mathcal{X}=\mathrm{Cliff}(f)=\K1\oplus\mathbb{V}$ be the Clifford Jordan algebra of $f$ over $\K$, see Definition \ref{def3}. Then $\mathcal{X}$ has a unique normalized trace $T$ such that $\mathcal{X}_0=\mathbb{V}$ and the product $*$ on $\mathcal{X}_0$ is 0. Now put $D_{v,v^\prime}=-[L_v,L_{v^\prime}]$, for $v,v^\prime\in\mathbb{V}$. {Here $L_v$ denotes the left multiplication by $v$}. Then one may identify $\fg=D_{\mathbb{V},\mathbb{V}}$, {where $D_{\mathbb{V},\mathbb{V}}$ is the $\bbbk$-span of $D_{v,v'}$, $v,v'\in\mathbb{V}$.}

Next assume that $\bb=\aa_{[n]}$ is the ring of Laurent polynomials in $n$ variables and $\mathcal{W}=\aa_{[n]}^{(m-1)}$ is the free $\aa_{[n]}$-module with basis $w_2,\ldots,w_m$. Fix $\tau_1,\ldots,\tau_m\in{\Lam:=\bbbz^n}$, $m\geq1$, such that 
$$\tau_1=0\quad\text{and}\quad\tau_i\not\equiv\tau_j (\mod\;2\bbbz^n)\quad\text{for}\quad i\neq j.$$ 
Define $g:\mathcal{W}\times\mathcal{W}\rightarrow\bb$ by
\begin{equation}\label{eq12}
	g\left(\sum_{i=2}^{m}a_iw_i,\sum_{i=2}^{m}b_iw_i\right)=\sum_{i=2}^{m}a_ib_ix^{\tau_i}.
\end{equation}
Let $\mathcal{Y}=\mathrm{Cliff}(g)=\bb1\oplus\mathcal{W}$ be the Clifford Jordan algebra of $g$ over $\K$, see Definition \ref{def3}. Then $\mathcal{Y}$ has a unique normalized trace $T$ such that $\mathcal{Y}_0=\mathcal{W}$ and the product $*$ on $\mathcal{Y}_0$ is 0. Set $D_{w,w^\prime}=-[L_w,L_{w^\prime}]$, for $w,w^\prime\in\mathcal{W}$. 

Now consider 
$$\ll:=\mathcal{T}(\mathrm{Cliff}(f)/\K,\mathrm{Cliff}(g)/\bb)=(\fg\otimes \bb)\oplus(\mathbb{V}\otimes\mathcal{W})\oplus D_{\mathcal{W},\mathcal{W}}.$$
Then $\ll$ has the root space decomposition $\ll=\sum_{\a\in\Delta}\ll_\a$ with respect to $\fh$, where
\begin{equation*}
	\ll_{\a}=\left\{\begin{array}{ll}
		\fg_\a\otimes \bb&\hbox{if }\a\in\Delta_{lg}\\
		\fg_\a\otimes \bb+\mathbb{V}_\a\otimes\mathcal{W}&\hbox{if }\a\in\Delta_{sh}\\
		\fh\otimes \bb+\mathbb{V}_0\otimes\mathcal{W}+D_{\mathcal{W},\mathcal{W}}&\hbox{if }\a=0.
	\end{array}\right.
\end{equation*}
Also $\mathcal{Y}=\mathrm{Cliff}(g)$ has the following $\Lam$-grading as a Jordan algebra:
\begin{equation*}
	\deg(x^\lam1)=2\lam
	\quadd\text{and}\quadd
	\deg(x^\lam w_i)=2\lam+\tau_i,
\end{equation*}
for $\lam\in\Lam$ and $2\leq i\leq m$ (see \cite[Proposition 5.26]{AG01}). This gives $\ll$ a $\Lam$-grading as a Lie algebra by
\begin{equation}\label{eq13}
	\ll^\lam=(\fg\otimes \bb^\lam)\oplus(\mathbb{V}\otimes\mathcal{W}^\lam)\oplus D_{\mathcal{W},\mathcal{W}}^\lam
\end{equation}
for $\lam\in\Lam$, where $\bb^\lam=\mathcal{Y}^\lam\cap\bb$, $\mathcal{W}^\lam=\mathcal{Y}^\lam\cap \mathcal{W}$ and $D_{\mathcal{W},\mathcal{W}}^\lam=\sum_{\mu+\nu=\lam}[L_{\mathcal{W}^\mu},L_{\mathcal{W}^\nu}]$. Then $\ll$ is $Q\times\Lam$-graded. We are now ready to state the coordinatization theorem for type $B_\ell (\ell\geq3)$ which is due to B. Allison and Y. Gao.

\begin{thm}\cite[Theorem 5.29]{AG01}\label{thm6}
	Any centerless Lie $\Lam$-torus of type $B_\ell, \ell\geq3$ is graded-isomorphic to $\mathcal{T}(\mathrm{Cliff}(f)/\K,\mathrm{Cliff}(g)/\aa_{[n]})$ where $f$ is a nondegenerate symmetric bilinear form on a	$(2\ell+1)$-dimensional vector space $\mathbb{V}$ over $\K$ and $g:\aa_{[n]}^{(m-1)}\times \aa_{[n]}^{(m-1)}\rightarrow \aa_{[n]}$ is defined by (\ref{eq12}).
\end{thm}

\begin{pro}\label{pro7}
	Assume that $\ll$ is a centerless $\Lam$-torus of type $B_\ell, \ell\geq3$. Then $\ll$ admits a Chevalley involution.	
\end{pro}
\proof
By Theorem \ref{thm6}, it is enough to show the existence of a Chevalley involution for $$\ll=\mathcal{T}(\mathrm{Cliff}(f)/\K,\mathrm{Cliff}(g)/\aa_{[n]})=(\fg\otimes\aa_{[n]})\oplus(\mathbb{V}\otimes\mathcal{W})\oplus D_{\mathcal{W},\mathcal{W}}$$
where $\mathcal{W}=\aa_{[n]}^{(m-1)}$ and $f:\mathbb{V}\times\mathbb{V}\rightarrow\K$, $g:\mathcal{W}\times\mathcal{W}\rightarrow\aa_{[n]}$ are as in Theorem \ref{thm6}.

Let $\tau:\mathrm{Cliff}(g)\rightarrow\mathrm{Cliff}(g)$ be the linear map given by
$$x^\lam1\mapsto x^{-\lam}1,\quadd x^\lam w_i\mapsto x^{-\lam-\tau_i} w_i$$
for $\lam\in\bbbz^n$ and $2\leq i\leq m$. One can easily see that $\tau$ is a pre-Chevalley involution for $\mathrm{Cliff}(g)$.

Next define $\bar\tau:\ll\rightarrow\ll$ by
$$\bar\tau(y\otimes x^\lam+v\otimes w+D_{w^\prime,w^{\prime\prime}})=\theta(y)\otimes x^{-\lam}+\overline{v}\otimes\tau(w)+D_{\tau(w^\prime),\tau(w^{\prime\prime})}$$
for $y\in\fg,\lam\in\Lam,v\in\mathbb{V},w,w^\prime,w^{\prime\prime}\in\mathcal{W}$, where $\bar{}:\mathbb{V}\rightarrow\mathbb{V}$ is the period $2$ linear map defined by 
$$\overline{v_i}=v_{l+i},\;\; i=1,\ldots,\ell,\quadd\overline{v_{2\ell+1}}=v_{2\ell+1},$$
and $\theta$ is a Chevalley involution of $\fg$. It is straightforward to check that $\bar\tau$ is a Lie algebra involution of $\ll$. Also considering the $\Lam$-grading (\ref{eq13}) of $\ll$,  we clearly have $\bar\tau(\ll^\lam)=\ll^{-\lam}$, for $\lam\in\Lam$. So $\bar\tau$ is a pre-Chevalley involution on $\ll$.

Finally, note that $\ll_0^0=\fh\otimes1\oplus D_{\mathcal{W},\mathcal{W}}^0$ with $D_{\mathcal{W},\mathcal{W}}^0=\sum_{\mu\in\Lam}[L_{\mathcal{W}^\mu},L_{\mathcal{W}^{-\mu}}]$. Then clearly $\bar\tau_{|_{\ll_0^0}}=-\id$. Thus $\bar\tau$ is a Chevalley involution of $\ll$ and this completes our proof.
\qed

\subsubsection{\bf{Type} $\mathbf{F_4}$}
Consider the root system $\Delta=\{0\}\cup\Delta_{sh}\cup\Delta_{lg}$
of type $F_4$ with $\Delta_{sh}$ and $\Delta_{lg}$ as short and long roots, respectively, and let $Q=Q(\Delta)$ be the root lattice of $\Delta$.

Let $\mathbb{J}$ be the {\it Albert algebra} over $\K$, the $27$-dimensional exceptional simple Jordan algebra over $\K$, and let $\fg=\mathrm{Der}(\mathbb{J})$. Then $\fg$ is known to be the finite dimensional simple Lie algebra of type $F_4$. Now let $T$ be the usual normalized trace on $\mathbb{J}$ and, for $x,x^\prime\in\mathbb{J}_0$, let $D_{x,x^\prime}=[L_x,L_{x'}]$. Then we have $\fg=D_{\mathbb{J}_0,\mathbb{J}_0}$ (see \cite{Jac71} and \cite[Chapter IX]{Jac68} for these facts).

Next we need to recall from \cite{BGKN95}  the definition of some alternative algebras (also see \cite[Example III.5.5]{AABGP97}). Assume that $\aa_0=\aa_{[n]}$ is the ring of Laurent polynomials in $n$ variables. Moreover, for $1\leq p\leq3$ with $p\leq n$, let
$$\aa_p=\aa_{p-1}\oplus t_p\aa_{p-1}$$
be the algebra obtained from $\aa_{p-1}$ using the Cayley-Dickson process with $t^2_p=x_p$. We also let $T$ be the normalized trace on $\aa_p$. Note that the last algebra $\aa_3$ is the octonion torus, see Example \ref{ex1}.

Next consider 
$$\ll:=\mathcal{T}(\mathbb{J}/\K,\mathcal{C}/\bb)=(\fg\otimes \bb)\oplus(\mathbb{J}_0\otimes\mathcal{C}_0)\oplus D_{\mathcal{C}_0,\mathcal{C}_0}$$
where $\bb=\aa_{[n]}$, $\cc=\aa_p$ for some $0\leq p\leq3$ and $D_{\mathcal{C}_0,\mathcal{C}_0}$ is the $\bbbk$-span of inner derivations
$$D_{y,y'}=\frac{1}{4}(L_{[y,y']}-R_{[y,y']}-3[L_y,R_{y'}])$$
of $\cc$, for $y,y'\in\cc_0$. Let $\Lam=\Z^n=\Z\lam_1\oplus\cdots\oplus\Z\lam_n$ where $\lam_1,\ldots,\lam_n$ is the standard basis of $\Lam$. Then $\ll$ is a $Q\times\Lam$-graded Lie algebra as follows:
Let $\fh$ be a Cartan subalgebra of $\fg$. First note that $\mathbb{J}_0$ is an irreducible $\fg$-module \cite[p.98]{Jac71} with the set of weights $\{0\}\cup\Delta_{sh}$ and
\begin{equation*}
	\dim_\K(\mathbb{J}_0)_\a=\left\{\begin{array}{ll}
		1&\hbox{if }\a\in\Delta_{sh}\\
		2&\hbox{if }\a=0
	\end{array}\right.
\end{equation*}
where $(\mathbb{J}_0)_\a$ is the the weight space of $\mathbb{J}_0$ corresponding to $\a$ \cite[$\S5, \S6$]{Jac71}. Then $\ll$ has the root space decomposition $\ll=\sum_{\a\in\Delta}\ll_\a$ with respect to $\fh$, where
\begin{equation*}
	\ll_{\a}=\left\{\begin{array}{ll}
		\fg_\a\otimes \bb&\hbox{if }\a\in\Delta_{lg}\\
		\fg_\a\otimes \bb+(\mathbb{J}_0)_\a\otimes\cc_0&\hbox{if }\a\in\Delta_{sh}\\
		\fh\otimes \bb+(\mathbb{J}_0)_0\otimes\cc_0+D_{\cc_0,\cc_0}&\hbox{if }\a=0.
	\end{array}\right.
\end{equation*}
The algebra $\cc$ has a unique $\Lam$-grading such that
\begin{equation*}
	\deg(t_i)=\lam_i,1\leq i\leq p,
	\quadd\text{and}\quadd
	\deg(x_i)=\lam_i,p+1\leq i\leq n,
\end{equation*}
see \cite[Proposition 5.26]{AG01}. This makes $\ll$ a $\Lam$-graded Lie algebra with
\begin{equation}\label{eq14}
	\ll^\lam=(\fg\otimes \bb^\lam)\oplus(\mathbb{J}_0\otimes\cc_0^\lam)\oplus D_{\cc_0,\cc_0}^\lam
\end{equation}
for $\lam\in\Lam$, where $\bb^\lam=\cc^\lam\cap \bb$, $\cc_0^\lam=\cc^\lam\cap\cc_0$ and $D_{\cc_0,\cc_0}^\lam=\sum_{\mu+\nu=\lam}D_{\cc_0^\mu,\cc_0^\nu}$. We are now ready to state the coordinatization theorem for type $F_4$.

\begin{thm}\cite[Theorem 5.50]{AG01}\label{thm7}
	Any centerless Lie $\Lam$-torus of type $F_4$ is graded-isomorphic to $\mathcal{T}(\mathbb{J}/\K,\mathcal{C}/\aa_{[n]})$ where $\mathbb{J}$ is the Albert algebra over $\K$ and $\cc=\aa_p$ for some $0\leq p\leq3$.
\end{thm}

\begin{pro}\label{pro8}
	Assume that $\ll$ is a centerless $\Lam$-torus of type $F_4$. Then $\ll$ admits a Chevalley involution.	
\end{pro}
\proof
{By Theorem \ref{thm7}, we may assume that} $$\ll=\mathcal{T}(\mathbb{J}/\K,\mathcal{C}/\aa_{[n]})=(\fg\otimes \aa_{[n]})\oplus(\mathbb{J}_0\otimes\mathcal{C}_0)\oplus D_{\mathcal{C}_0,\mathcal{C}_0}$$
where $\mathbb{J}$ is the Albert algebra over $\K$ and $\cc=\aa_p$ for some $0\leq p\leq3$.

Let $\tau:\aa_p\rightarrow\aa_p$ be the {linear} map defined by
\begin{equation*}
	t_i\mapsto t_i^{-1},1\leq i\leq p,
	\quadd\text{and}\quadd
	x_i\mapsto x_i^{-1},p+1\leq i\leq n.
\end{equation*}
This together with the way the $\Lam$-grading on $\aa_p$ is defined implies that $\tau$ is a pre-Chevalley involution of $\aa_p$ (for $\aa_3$ see Proposition \ref{pro2}). 

Next fix a basis $\{a_1,a_2,a_\a\mid a_1,a_2\in(\mathbb{J}_0)_0,a_\a\in(\mathbb{J}_0)_\a,\a\in\Delta_{sh}\}$ for $\mathbb{J}_0$ and let $\bar{}:\mathbb{J}_0\rightarrow\mathbb{J}_0$ be the period $2$ linear map such that
$$\overline{a_i}=a_i,i=1,2,\quadd\overline{a_\a}=a_{-\a}.$$
Then we define $\bar\tau:\ll\rightarrow\ll$ by
$$\bar\tau(y\otimes x+a\otimes c+D_{c^\prime,c^{\prime\prime}})=\theta(y)\otimes \tau(x)+\overline{a}\otimes\tau(c)+D_{\tau(c^\prime),\tau(c^{\prime\prime})}$$
for $y\in\fg,x\in\aa_{[n]},a\in\mathbb{J}_0,c,c^\prime,c^{\prime\prime}\in\cc_0$, where $\theta$ is a Chevalley involution of $\fg$. One directly checks that $\bar\tau$ is a Lie algebra involution of $\ll$. Also considering the $\Lam$-grading (\ref{eq14}) of $\ll$,  we clearly have $\bar\tau(\ll^\lam)=\ll^{-\lam}$, for $\lam\in\Lam$. So $\bar\tau$ is a pre-Chevalley involution on $\ll$.

Finally, note that $\ll_0^0=\fh\otimes1\oplus D_{\cc_0,\cc_0}^0$. But $D_{\cc_0,\cc_0}^0=\{0\}$, see \cite[Proposition III.5.35 (d)]{AABGP97}. Then clearly $\bar\tau_{|_{\ll_0^0}}=-\id$. Thus $\bar\tau$ is a Chevalley involution of $\ll$ and this completes our proof.
\qed

\subsubsection{\bf{Type} $\mathbf{G_2}$}

Assume that $\Delta=\{0\}\cup\Delta_{sh}\cup\Delta_{lg}$
is the finite root system of type $G_2$ with $\Delta_{sh}$ and $\Delta_{lg}$ as short and long roots, respectively and let $Q=Q(\Delta)$ be the root lattice of $\Delta$.

Let $\mathbb{A}$ be the {\it Cayley algebra} over $\K$, the $8$-dimensional simple nonassociative alternative algebra over $\K$, and let $\fg=\mathrm{Der}(\mathbb{A})$. Then $\fg$ is known to be the finite dimensional simple Lie algebra of type $G_2$. Now let $T$ be the usual normalized trace on $\mathbb{A}$ and, for $x,x^\prime\in\mathbb{J}_0$, let 
$$D_{x,x^\prime}=\frac{1}{4}(L_{[x,x']}-R_{[x,x']}-3[L_x,R_{x'}]).$$ 
Then we have $\fg=D_{\mathbb{A}_0,\mathbb{A}_0}$ (see \cite{Jac71} and \cite[Chapter III]{Sch66} for these facts).

Next we recall the definition of the Jordan algebras $\jj_p$ for $0\leq p\leq3$ with $p\leq n$ (see \cite[Example III.5.10]{AABGP97}). Let $\jj_0=\aa_{[n]}$
{be the algebra of Laurent polynomials  in variables $x_1,\ldots, x_n$.} Let
$\jj_1$ be the commutative associative algebra over $\aa_{[n]}$ with generator $t_1$ subject to the relation $t_1^3=x_1$, and $\jj_2$ be the plus algebra of the associative algebra over $\aa_{[n]}$ generated by $t_1,t_2$ subject to the relations $t_1^3=x_1,t_2^3=x_2,t_1t_2=e^{2\pi i/3}t_2t_1$. Now let $$\jj_3=\jj_2\oplus\jj_2\cdot t_3\oplus\jj_2\cdot t_3^2$$ 
be the Jordan algebra obtained from $\jj_2$ using Tits first Jordan algebra construction \cite[Chapter IX]{Jac68} with $t_3^3=x_3$. We also let $T$ be the normalized trace on $\jj_p$. We mention that the last algebra $\jj_3$ is the Albert torus, see Example \ref{ex4}(iv).

Next consider 
$$\ll:=\mathcal{T}(\mathbb{A}/\K,\mathcal{J}/\bb)=(\fg\otimes \bb)\oplus(\mathbb{A}_0\otimes\mathcal{J}_0)\oplus D_{\mathcal{J}_0,\mathcal{J}_0}$$
where $\bb=\aa_{[n]}$, $\jj=\jj_p$ for some $0\leq p\leq3$ and $D_{\mathcal{J}_0,\mathcal{J}_0}$ is the $\bbbk$-span of inner derivations $D_{y,y'}=[L_y,L_{y'}]$ of $\jj$, for $y,y'\in\jj_0$. 
Let $\Lam=\Z^n=\Z\lam_1\oplus\cdots\oplus\Z\lam_n$ where $\lam_1,\ldots,\lam_n$ is the standard basis of $\bbbz^n$. Then $\ll$ is a $Q\times\Lam$-graded Lie algebra as follows:
Let $\fh$ be a Cartan subalgebra of $\fg$. First note that $\mathbb{A}_0$ is an irreducible $\fg$-module \cite[Proposition 4]{Jac71} with the set of weights $\{0\}\cup\Delta_{sh}$ and $\dim_\K(\mathbb{A}_0)_\a=1$ for all $\a\in\{0\}\cup\Delta_{sh}$, where $(\mathbb{A}_0)_\a$ is the the weight space of $\mathbb{A}_0$ corresponding to $\a$ \cite[$\S2$]{Jac71}. Then $\ll$ has the root space decomposition $\ll=\sum_{\a\in\Delta}\ll_\a$ with respect to $\fh$, where
\begin{equation*}
	\ll_{\a}=\left\{\begin{array}{ll}
		\fg_\a\otimes \bb&\hbox{if }\a\in\Delta_{lg}\\
		\fg_\a\otimes \bb+(\mathbb{A}_0)_\a\otimes\jj_0&\hbox{if }\a\in\Delta_{sh}\\
		\fh\otimes \bb+(\mathbb{A}_0)_0\otimes\jj_0+D_{\jj_0,\jj_0}&\hbox{if }\a=0.
	\end{array}\right.
\end{equation*}
The algebra $\jj$ has a unique $\Lam$-grading such that
\begin{equation*}
	\deg(t_i)=\lam_i,1\leq i\leq p,
	\quadd\text{and}\quadd
	\deg(x_i)=\lam_i,p+1\leq i\leq n,
\end{equation*}
see \cite[Proposition 5.58]{AG01}. This gives $\ll$ a $\Lam$-grading as a Lie algebra with
\begin{equation}\label{eq15}
	\ll^\lam=(\fg\otimes \bb^\lam)\oplus(\mathbb{A}_0\otimes\jj_0^\lam)\oplus D_{\jj_0,\jj_0}^\lam
\end{equation}
for $\lam\in\Lam$, where $\bb^\lam=\jj^\lam\cap \bb$, $\jj_0^\lam=\jj^\lam\cap\jj_0$ and $D_{\jj_0,\jj_0}^\lam=\sum_{\mu+\nu=\lam}[L_{\jj_0^\mu},L_{\jj_0^\nu}]$. Then $\ll$ is $Q\times\Lam$-graded. We are now ready to state the coordinatization theorem for type $G_2$ which is due to B. Allison and Y. Gao.

\begin{thm}\cite[Theorem 5.63]{AG01}\label{thm8}
	Any centerless Lie $\Lam$-torus of type $G_2$ is graded-isomorphic to $\mathcal{T}(\mathbb{A}/\K,\jj/\aa_{[n]})$ where $\mathbb{A}$ is the Cayley algebra over $\K$ and $\jj=\jj_p$ for some $0\leq p\leq3$.
\end{thm}

\begin{pro}\label{pro9}
	Assume that $\ll$ is a centerless $\Lam$-torus of type $G_2$. Then $\ll$ admits a Chevalley involution.	
\end{pro}
\proof
By Theorem \ref{thm8}, it is enough to show the existence of a Chevalley involution for $$\ll=\mathcal{T}(\mathbb{A}/\K,\jj/\aa_{[n]})=(\fg\otimes \aa_{[n]})\oplus(\mathbb{A}_0\otimes\jj_0)\oplus D_{\jj_0,\jj_0}$$
where $\mathbb{A}$ is the Cayley algebra over $\K$ and $\jj=\jj_p$ for some $0\leq p\leq3$.

Let $\tau:\jj_p\rightarrow\jj_p$ be the linear map defined by
\begin{equation*}
	t_i\mapsto t_i^{-1},1\leq i\leq p,
	\quadd\text{and}\quadd
	x_i\mapsto x_i^{-1},p+1\leq i\leq n.
\end{equation*}
Now one can easily see that $\tau$ is a pre-Chevalley involution for $\jj_p$. 

Next fix a basis $\{a_\a\mid a_\a\in(\mathbb{A}_0)_\a,\a\in\{0\}\cup\Delta_{sh}\}$ for $\mathbb{A}_0$ and let $\bar{}:\mathbb{A}_0\rightarrow\mathbb{A}_0$ be the period $2$ linear map such that $\overline{a_\a}=a_{-\a}$.
Then define $\bar\tau:\ll\rightarrow\ll$ by
$$\bar\tau(y\otimes x+a\otimes z+D_{z^\prime,z^{\prime\prime}})=\theta(y)\otimes \tau(x)+\overline{a}\otimes\tau(z)+D_{\tau(z^\prime),\tau(z^{\prime\prime})}$$
for $y\in\fg,x\in\aa_{[n]},a\in\mathbb{A}_0,z,z^\prime,z^{\prime\prime}\in\jj_0$, where $\theta$ is a Chevalley involution of $\fg$. It is straightforward to check that $\bar\tau$ is a Lie algebra involution of $\ll$. Also considering the $\Lam$-grading (\ref{eq15}) of $\ll$,  we clearly have $\bar\tau(\ll^\lam)=\ll^{-\lam}$, for $\lam\in\Lam$. So $\bar\tau$ is a pre-Chevalley involution on $\ll$.

Finally, note that $\ll_0^0=\fh\otimes1\oplus D_{\jj_0,\jj_0}^0$. But $D_{\jj_0,\jj_0}^0=\{0\}$, see \cite[Proposition III.5.22 (d)]{AABGP97}. Then clearly $\bar\tau_{|_{\ll_0^0}}=-\id$. Thus $\bar\tau$ is a Chevalley involution of $\ll$ and this completes our proof.
\qed

\section{\bf Chevalley involutions for fgc centerless tori}\setcounter{equation}{0}\label{sec7}
{In this section, we discuss the existence of Chevalley involutions for fgc centerless Lie tori. We recall from Definition \ref{fgc} that a centerless Lie torus which is finitely generated as a module over its centroid is called fgc. According to Theorem \ref{thm9}(b), all centerless Lie tori of types different from $A_\ell$ are fgc.
Moreover, according to Theorem \cite[Theorem 3.3.1]{ABFP09} (see Theorem \ref{fgc1} below), any fgc centerless Lie torus is a ''multi-loop algebra''.
Therefore, the existence of Chevalley involutions for centerless Lie tori amounts to the existence of Chevalley involutions for multi-loop algebras. } 

\subsection{Multi-loop Lie $\bbbz^n$-tori}\setcounter{equation}{0}
{Let $\fg$ be a finite dimensional simple Lie algebra and $\sg_1,\ldots,\sg_\nu$, be commuting finite order automorphisms of $\fg$ with periods
		$m_1,\ldots,m_\nu$ respectively. We set $\pmb\sg=(\sg_1,\ldots,\sg_\nu)$. 
		Let $\omega_i$ be a primitive $m_i^{th}$-root of
		unity for $1\leq i\leq\nu$. Then 
		$$
		\fg=\bigoplus_{\lam\in\bbbz^n}\fg^{\bar\lam},
		$$ where for $\lam=(\lam_1,\ldots,\lam_\nu)\in\Lam:=\bbbz^\nu$,
		$\bar\lam:=(\bar\lam_1,\ldots,\bar\lam_\nu)$ with 
		$\bar\lam_j:=\lam_j+m_j\bbbz\in\bbbz_{m_j}$, and
		$$
		\fg^{\bar\lam}=\{x\in\fg\mid\sg_j(x)=\omega^{\lam_j}_{j}x\hbox{
			for }1\leq j\leq\nu\}.
		$$
		Let $\pi_{\bar{\lam}}:\fg\rightarrow\fg^{\bar{\lam}}$ denote the projection onto $\fg^{\bar{\lam}}$.
		Let $\aa=\bbbk[z^{\pm1}_1,\ldots,z^{\pm1}_\nu]$ be the algebra of Laurent polynomials in $\nu$-variables equipped with the natural $\Lam$-grading $\aa=\sum_{\lam\in\Lam}\bbbk z^\lam$, $z^\lam=z_1^{\lam_1},\ldots, z_\nu^{\lam_\nu}$, $\lam=(\lam_1,\ldots,\lam_\nu).$}
		
		{(i) The subalgebra
		\begin{equation}\label{ref}
			M(\fg,\pmb\sg):=\bigoplus_{\lam\in\bbbz^n}\pi_{\bar{\lam}}(\fg)\otimes
			z^\lam=
			\bigoplus_{\lam\in\bbbz^n}\fg^{\bar\lam}\otimes
			z^\lam
		\end{equation}
		of $L(\fg,\aa):=\fg\otimes\aa$ is called the {\it $\nu$-step multi-loop
			algebra} based on $\pmb\sg$ and $\fg$. }
		
		{(ii) It is known that $\fg^{\pmb\sg}$ contains a subalgebra $\fh'$ which is $\fg$-ad-diagonalizable, see {\cite[Remark 3.2.2(ii)]{ABFP09}}. Let $\Delta'$ be the set of weights of $\fg$ with respect to $\fh'$. Set
		$$\hbox{LT}(\fg,\pmb\sg,\fh'):=\sum_{\lam\in\Lam,\a'\in\Delta'}\fg_{\a'}^{\bar\lam}\otimes z^\lam,$$ 
		where $\fg^{\bar\lam}_{\a'}$ is the intersection of $\fg^{\bar\lam}$ and the $\a'$-weight space of $\fg$, with respect to $\fh'$. Then $\hbox{LT}(\sg,\pmb\sg,\fh')$ is a subalgebra of $\fg\otimes\aa$
		with the compatible $(\la\Delta'\ra\times\Lam)$-grading given above.
		We note that as $\Lam$-graded Lie algebras, $\hbox{LT}(\fg,\pmb\sg,\fh')=M(\fg,\pmb\sg)$.}

	{\begin{thm}\label{fgc1}\cite[Theorem 3.3.1]{ABFP09}
		Let $E(\ll,D,\kappa)$ be a fgc extended affine Lie algebra of nullity $\nu\geq 1$, where
		$\ll$ is a centerless Lie torus. Then $\ll$ is bi-isomorphic to a multi-loop Lie algebra
		$LT(\fg,\pmb\sg,\fh')=M(\fg,\pmb\sg)$ for some $\fg,\pmb\sg,\fh'$, with
		$\ll^0_0\cong\fh'\otimes 1$.
				\end{thm}
}	
	
\subsection{Chevalley involutions for fgc multi-loop algebras}\setcounter{equation}{0}\label{subsec7}
{We discuss Chevalley involutions for  fgc centerless Lie tori.}
In what follows for $\pmb\sg=(\sg_1,\ldots,\sg_\nu)$ we write $\pmb\sg\tau=\tau\pmb\sg$ to indicate that an automorphism $\tau$ commutes with all $\sg_i$.
\begin{pro}\label{fgc9} 
{Let $\ll=\hbox{LT}(\fg,\pmb\sg,\fh')=
M(\fg,\pmb\sg)$ be a fgc centerless Lie torus. Let $\tau$ be a Chevalley involution for $\fg$ such that $\tau\pmb\sg=\pmb\sg\tau$ and $\tau(h')=-h'$ for $h'\in\fh'$. 
Assume that $\fg$ has an automorphism $\psi$  satisfying:
}

- $\psi(\fg^{\bar\lam})\sub\fg^{-\bar\lam}$, for each $\lam$,

- $\psi_{|_{\fh'}}=\id_{|_{\fh'}}$. \\
{Then the assignment $\bar\tau_\psi:
 x\otimes a^\lam\longmapsto
\psi\tau(x)\otimes a^{-\lam},$
$x\in\fg$, $\lam\in\Lam$
defines a pre-Chevalley involution  $\bar\tau_\psi$ of $\fg\otimes\aa$ which restricts to a Chevalley involution for
$\ll$.}
\end{pro}

\proof {Let $\psi$ and $\fh'$ be as in the statement. We have
$\ll^0_0=\fh'\otimes1.$
Clearly, $\bar\tau_\psi$ is a pre-Chevalley involution on $\fg\otimes\aa$.
Since $\tau$ commutes with $\pmb\sg$, we have $\tau(\pi_\lam(\fg))=\pi_\lam(\fg)$, for $\lam\in\Lam$. This together with  $\psi(\fg^\lam)\sub\fg^{-\lam}$, gives
$$\bar\tau_\psi(\pi_\lam(\fg)\otimes a^\lam)=\psi\tau(\pi_\lam(\fg))\otimes a^{-\lam}
=\psi(\pi_\lam(\fg))\otimes a^{-\lam}
\sub\pi_{-\lam}(\g)\otimes a^{-\lam}.$$
Thus $\bar\tau_\psi$ restricts to an automorphism of $\ll$, satisfying
$\bar\tau_\psi(\ll^\lam)\sub\ll^{-\lam}.$
Moreover,  for $h'\in\fh'$, 
$
\bar\tau_\psi(h'\otimes 1)=
\psi(\tau(h'))\otimes1=\psi(-h')\otimes1=-h'\otimes1.
$
Thus ${\bar{\tau}_\psi}$ is a Chevalley involution for $\ll$.\qed}

\begin{cor}\label{fgc15}
{Let $\ll=\hbox{LT}(\fg,\pmb\sg,\fh')=
M(\fg,\pmb\sg)$ be a fgc centerless Lie torus, where all $\sg_i$'s are graph automorphisms.
Assume that $\fg$ has an automorphism $\psi$  satisfying:}

- $\psi(\fg^{\bar\lam})\sub\fg^{-\bar\lam}$, for each $\lam$,

- $\psi_{|_{\fh'}}=\id_{|_{\fh'}}$. \\
{Then $\ll$ is equipped with a Chevalley involution. 
In particular, if each $\sg_i$ has order $1$, or $2$, then the statement holds with $\psi=\id$. }
\end{cor}

\proof 
{We consider a Cartan subalgebra $\fh$ of $\fg$ containing $\fh'$, and a Chevalley involution
$\tau$ for $(\fg,\fh)$. Since Chevalley involutions commute with graph automorphisms, we have
$\tau\pmb\sg=\pmb\sg\tau$. Thus the first part of the statement follows from Proposition \ref{fgc9}. If each $\sg_i$ has order $1$ or $2$, then $\bar\lam=-\bar\lam$,  for each $\lam$, and so the second part of the statement is clear.\qed}

\begin{rem}\label{fgc10}
{We discuss the existence of automorphism $\psi$ in Corollary \ref{fgc15}. As it is mentioned in this corollary,  if each $\sg_i$ has order $1$ or $2$, then $\psi=\id$ works.
This already covers many of the possible cases.
For graph automorphisms of order $3$ which accuses only in  type $D_4$, a construction of $\psi$ discussed in \cite[Corollary 3.29]{AFI22} should be considered. }
\end{rem}     

{The following is perhaps the easiest example to see how the ad-diagonalizable subalgebra $\fh'$ of
$\fg^\sg$ appearing in $\hbox{LT}(\fg,\pmb\sg,\fh')$ applies. }	  
\begin{exa}
{Let $\fg=\frak{sl}_2(\bbbk)$ with the standard Chevalley basis $\{e,h,f\}$.	
Consider the Chevalley involution $\sg:e\leftrightarrow -f$, $h\rightarrow -h$.
Then $\fg^{\bar0}=\fg^\sg=\bbbk(e-f)$ and $\fg^{\bar 1}=\bbbk h\oplus\bbbk(e+f)$.
Take $h'=i(e-f)/2$ and $\fh'=\bbbk h'$. Then $\ad h'$ is semi-simple on $\fg$. In fact 
if $y=e+f-ih$ and $z=e+f+ih$, then $[h',y]=-y$, $[h',z]=z$, so we have the eigenspace decomposition
$\fg=\fg_0+\fg_1+\fg_{-1}$ with $\fg_0=\fh'$, $\fg_1=\bbbk z$ and $\fg_{-1}=\bbbk y$. 	
Now we have $M(\fg,\sg)=(\fg^\sg\otimes \bbbk[z^{\pm 2}])\oplus
(\fg^{\bar{1}}\otimes z\bbbk[z^{\pm2}])=(\fh'\otimes\bbbk[z^{\pm2}])\oplus
((\bbbk y\oplus\bbbk z)\otimes z\bbbk[z^{\pm2}]).$ Next, we consider the Chevalley involution 
$\tau:y\leftrightarrow -z$, $h'\rightarrow -h'$ for $\fg$. We have $\tau\sg=\sg\tau$. Then
Proportion \ref{fgc9}, applies with $\psi=\id.$}
\end{exa}


\end{document}